\documentclass[12pt]{amsart}
\usepackage{latexsym,verbatim}
\usepackage{euscript}
\usepackage{amsfonts,amscd,amssymb,epsfig,latexsym,amsmath}
\usepackage{graphicx}

\newenvironment{pf}{\begin{trivlist} \item[] {\it Proof.} \ }{\hfill
\qed \end{trivlist} }

\newtheorem{theorem}{Theorem}[section]
\newtheorem{lemma}[theorem]{Lemma}
\newtheorem{proposition}[theorem]{Proposition}
\newtheorem{corollary}[theorem]{Corollary}

\newtheorem{propvoid}{Proposition} 
\newtheorem{definition}[theorem]{Definition}
\newtheorem{ex}[theorem]{Example}
\newenvironment{example}{\begin{ex} \rm}{\end{ex}}

\renewcommand{\caption}{}

\numberwithin{equation}{section}

\def\&{\wedge}
\newcommand{\w}{\omega}
\newcommand{\alp}{\alpha}

\newcommand{\eps}{\varepsilon}

\newcommand{\calI}{{\mathcal I}}
\newcommand{\calJ}{{\mathcal J}}

\newcommand{\calL}{{\mathcal L}}
\newcommand{\calQ}{{\mathcal Q}}

\newcommand{\scB}{{\EuScript B}}

\newcommand{\scH}{{\EuScript H}}

\newcommand{\scL}{{\EuScript L}}
\newcommand{\scU}{{\EuScript U}}
\newcommand{\scX}{{\EuScript X}}
\newcommand{\scY}{{\EuScript Y}}
\newcommand{\scZ}{{\EuScript Z}}
\newcommand{\scM}{{\EuScript M}}
\newcommand{\hook}{\, {\rule[0in]{2mm}{0.25mm}\rule{0.25mm}{2mm}} \, }
\newcommand{\bb}{\mathbb}
\newcommand{\mbf}{\mathbf}

\begin{document}

\title{Sub-Finsler geometry in dimension three}

\author{Jeanne N. Clelland}
\address{Department of Mathematics, 395 UCB, University of
Colorado\\
Boulder, CO 80309-0395}
\email{Jeanne.Clelland@colorado.edu}
\author{Christopher G. Moseley}
\address{Department of Mathematical Sciences\\ U.S. Military Academy\\ West
Point, NY 10996}
\email{Chris.Moseley@usma.edu}


\begin{abstract}

We define the notion of {\em sub-Finsler geometry} as a natural 
generalization of sub-Riemannian geometry with applications to 
optimal control theory.  We compute a complete set of local 
invariants, geodesic equations, and the Jacobi operator for the 
three-dimensional case and investigate homogeneous examples.
\end{abstract}

\maketitle

\section{Introduction}

Much attention has been given in recent years to sub-Riemannian
geometry; it is a rich subject with many applications.  In this paper
we introduce the notion of {\em sub-Finsler geometry}, a natural
generalization of sub-Riemannian geometry.

The motivation for this generalization comes from optimal control
theory.  A control system is usually presented in local coordinates
as an underdetermined system of ordinary differential equations
\begin{equation}
   \dot{x} = f(x, u), \label{gencontrol}
\end{equation}
where $x \in \bb{R}^n$ represents the \emph{state} of the system and $u \in
\bb{R}^s$ represents the {\em controls}, i.e., variables which may be
specified freely in order to ``steer" the system in a desired
direction.  More generally, $x$ and $u$ may take values in an
$n$-dimensional manifold $\scX$ and an $s$-dimensional manifold $\scU$,
respectively.  Typically there are constraints on
how the system may be ``steered" from one state to another, so that 
$s < n$.  The
systems of greatest interest are {\em controllable}, i.e., given any
two states $x_1, x_2$, there exists a solution curve of
\eqref{gencontrol} connecting $x_1$ to $x_2$.

Consider the large class of systems which are linear (but not
affine linear) in the control variables $u$ and depend smoothly on
the state variables $x$, i.e., systems of the form
\begin{equation}
   \dot{x} = f(x) u, \label{qlincontrol}
\end{equation}
where $f(x)$ is a matrix whose entries are arbitrary smooth functions of $x$.
This class is by no means all-inclusive, but it does contain many
systems of interest; an example is given below.  For such a system,
\emph{admissible} paths in the state space are those for which the tangent
vector to the path at each point $x \in \scX$ is contained in the
subspace $D_x \subset T_x\scX$ determined by the image of the $n \times
s$ matrix $f(x)$.  Often this matrix is smooth and has constant
rank $s$, in which case $D$ is a rank $s$ distribution on $\scX$.  (In
this case the variables $(x, u)$ may be regarded as local coordinates
on the distribution $(\scX, D)$.)  Thus the admissible paths in the
state space are precisely the {\em horizontal curves} of the
distribution $D$, i.e., curves whose tangent vectors at each point
are contained in $D$.  By a theorem of Chow \cite{Chow39}, the system
\eqref{qlincontrol} is controllable if and only if the distribution
$D$ on $\scX$ is {\em bracket-generating}, i.e., if the iterated
brackets of vector fields contained in $D$ span the entire tangent
space at each point $x \in \scX$.

Given a distribution $(\scX, D)$ representing a system of the form
\eqref{qlincontrol}, we next consider the problem of {\em optimal
control}: what is the most efficient path between two given points in
$\scX$?  In order to answer this question, we must have some measure of
the cost required to move in the state space.  This measure is
typically specified by a first-order Lagrangian functional $L$
defined on the horizontal curves of $D$: given a horizontal curve
$\gamma: [a, b] \to \scX$, the {\em action} $\scL(\gamma)$ is defined to be
\[ \scL(\gamma) = \int_{\gamma} L(x, \dot{x})\, dx = \int_{\gamma}
\bar{L}(x, u) \, dx, \]
where, since $\gamma$ is a solution curve of \eqref{qlincontrol}, we
define $\bar{L}(x, u) = L(x, f(x)u)$.  Often the Lagrangian has the
form
\[ \bar{L}(x, u) = \sqrt{g_{ij}(x) u^i u^j} \]
(summation on repeated indices being understood),
and in this case it defines a {\em sub-Rieman\-nian metric} $\langle,
\rangle$ on $D$ (i.e., a Riemannian metric on each subspace $D_x
\subset T_x\scX$) in the obvious way.  Horizontal paths which minimize
the action functional are  precisely the geodesics of the
sub-Riemannian metric.

\begin{example}\label{planewheel}
Consider a wheel rolling without slipping on the Euclidean plane
$\bb{E}^2$.  The wheel's configuration can be represented by the
vector ${}^t\hskip-1pt(x, y, \varphi, \psi)$, where $(x,y)$ is the
wheel's point of contact with the plane, $\phi$ is the angle of
rotation of a marked point on the wheel from the vertical, and $\psi$
is the wheel's heading angle, i.e., the angle made by the tangent
line to the curve traced by the wheel on the plane with the $x$-axis.
Thus the state space has dimension four and is naturally isomorphic to
$\bb{R}^2 \times S^1 \times S^1$.

The condition that the wheel rolls without slipping is equivalent to
the statement that its path ${}^t\hskip-1pt(x(t), y(t), \varphi(t),
\psi(t))$ in the state space satisfies the differential equation
\[ \begin{bmatrix} \dot{x} \\ \dot{y} \\ \dot{\varphi} \\ \dot{\psi}
\end{bmatrix} =
u_1(t) \begin{bmatrix} \cos \psi \\ \sin \psi \\ 1 \\ 0 \end{bmatrix} +
u_2(t) \begin{bmatrix} 0 \\ 0 \\ 0 \\ 1 \end{bmatrix} \]
for some control functions $u_1(t), u_2(t)$.  Thus the velocity
vector ${}^t\hskip-1pt(\dot{x}, \dot{y}, \dot{\varphi}, \dot{\psi})$
of any solution curve must lie in the distribution $D$ spanned by the
vector fields
\begin{align*}
V_1 & = (\cos \psi) \frac{\partial}{\partial x} + (\sin \psi)
\frac{\partial}{\partial y} + \frac{\partial}{\partial \varphi} \\
V_2 & = \frac{\partial}{\partial \psi}.
\end{align*}
A natural sub-Riemannian metric on $D$ is obtained by declaring the
vector fields $V_1,\, V_2$ to be orthonormal, i.e., by setting
\[ \langle u_1 V_1 + u_2 V_2,\, u_1 V_1 + u_2 V_2 \rangle  = u_1^2 + u_2^2. \]
The integral of this quadratic form measures the work done in rotating the heading
angle $\psi$ at the rate $\dot{\psi}$ and propelling the wheel forward at the rate $\dot{\varphi}$.
\end{example}

But what if the natural measure on horizontal curves is not the
square root of a quadratic form?  For instance, suppose we modified
Example \ref{planewheel} by rolling the wheel on an {\em inclined}
plane?  (Assume that the wheel has sufficient friction to remain
motionless if no energy is put into the system.)  We would expect
more energy to be required to move the wheel uphill than downhill, so
the natural Lagrangian would not even be symmetric in $u$
(i.e., it would not satisfy the condition $\bar{L}(x, -u) = \bar{L}(x, u)$), let
alone be the square root of a quadratic form in $u$.  It is not
difficult to imagine examples where the dependence of $\bar{L}$ on $u$ becomes quite complicated as $u$ changes direction.  This leads us to
generalize the notion of a sub-Riemannian metric on $(\scX, D)$ by
replacing the Riemannian metric on each subspace $D_x \subset T_x\scX$
with a Finsler metric.

Recall that a {\em Finsler metric} on a manifold $\scM$ is a function
\[ F: T\scM \to [0, \infty) \]
with the following properties:
\begin{enumerate}
\item{Regularity: $F$ is $C^{\infty}$ on the slit tangent bundle $T\scM
\setminus 0$.}
\item{Positive homogeneity: $F(x, \lambda y) = \lambda F(x, y)$ for
all $\lambda > 0$. (Here $x$ is any system of local coordinates on
$\scM$ and $(x,y)$ is the corresponding canonical coordinate system on
$T\scM$.)}
\item{Strong convexity: The $n \times n$ Hessian matrix
\[ \left[ \frac{\partial^2 (\tfrac{1}{2} F^2)}{\partial y^i\,
\partial y^j} \right] \]
is positive definite at every point of $T\scM \setminus 0$.}\label{strongcon}
\end{enumerate}
(For details, see \cite{BCS00}.)  In other words, a Finsler metric on
a manifold $\scM$ is a smoothly varying Minkowski norm on each tangent
space $T_x\scM$. 

Condition \ref{strongcon} implies that the ``unit
sphere" in each tangent space $T_x\scM$ (also known as the {\em indicatrix} for the Finsler metric on $T_x\scM$) is a smooth, strictly convex
hypersurface enclosing the origin $0_x \in T_x\scM$.  The converse is almost -- but not quite -- true: there exist strictly convex hypersurfaces for which the corresponding Hessian matrix is only positive semi-definite along a closed subset; see \cite{BCS00} for examples.  We will say that a hypersurface $\Sigma_x \subset T_x\scM$ which encloses the origin is {\em strongly convex} if it is the indicatrix for a Minkowski norm on $T_x\scM$; thus strong convexity implies strict convexity, but not vice-versa.
  
In the Riemannian case, the indicatrix must be an ellipsoid centered at $0_x$, but in the Finsler case it may be much more general.  In particular, it need not be symmetric about the origin.

We are now ready to define our primary object of study.

\begin{definition}
A {\em sub-Finsler metric} on a smooth distribution $D$ of rank $s$
on an $n$-dimensional manifold $\scX$ is a smoothly varying Finsler
metric on each subspace $D_x \subset T_x\scX$.  A {\em sub-Finsler
manifold}, denoted by the triple $(\scX, D, F)$, is a smooth
$n$-dimensional manifold $\scX$ equipped with a sub-Finsler metric $F$
on a bracket-generating distribution $D$ of rank $s>0$.  The {\em
length} of a horizontal curve $\gamma:[a, b] \to \scX$ is
\[ \calL(\gamma) = \int_a^b F(\dot{\gamma}(t))\, dt. \]
\end{definition}

Replacing the Riemannian metric on $D$ by a Finsler metric allows
more general action functionals to be considered.  The rather
stringent requirement that the Lagrangian be the square root of a
quadratic form is replaced by the more natural requirement that it
be positive-homogeneous in $u$ (which is necessary if the length of
an oriented curve is to be independent of parametrization), and that it be
strongly convex (which is necessary if there are to exist locally
minimizing paths in every direction).  The problem of finding
minimizing paths satisfying \eqref{qlincontrol} is equivalent to finding geodesics of the sub-Finsler manifold $(\scX, D,
F)$.

In this paper we will investigate sub-Finsler manifolds in the
simplest nontrivial case: a three-dimensional manifold $\scX$ with a 
rank two contact
distribution $D$.  We will work locally, and thus we will not generally concern
ourselves with the issue of local vs. global existence of objects
such as coordinates, vector fields, and differential forms.

In the next two sections we will review some results of Hughen
\cite{Hughen95} concerning sub-Riemannian geometry in dimension three and
some results of Cartan \cite{Bryant96, Cartan34} concerning the geometry of
Finsler surfaces.  We will then combine these techniques to construct a
complete set of local invariants for sub-Finsler manifolds in dimension
three via \'{E}lie Cartan's method of equivalence.  (See \cite{Gardner89} for
an exposition of this method.  The reader should be aware that where
Gardner uses left group actions, we use right group actions for
greater ease of computation.)  Additionally, we will derive the 
geodesic equations, compute the Jacobi operator for the second 
variation problem, and investigate homogeneous examples.

\section{Review of sub-Riemannian geometry of 3-manifolds}\label{subRiemsect}

The material in this section is taken from Keener Hughen's Ph.D.
thesis \cite{Hughen95}.  Unfortunately this thesis was never
published, but some of the results are summarized in \cite{Montgomery02}.

Let $(\scX, D, \langle, \rangle)$ be a sub-Riemannian structure on a
$3$-manifold $\scX$ with a contact distribution $D$.  A local coframing
$(\eta^1, \eta^2, \eta^3)$ on $\scX$ is said to be {\em $0$-adapted} to
the sub-Riemannian structure if $D = \{\eta^3\}^{\perp}$ and
$\langle, \rangle = (\eta^1)^2 + (\eta^2)^2$.  The set of $0$-adapted
coframings of $\scX$ forms a $G_0$-structure $\scB_0 \to \scX$, where
$G_0$ is the Lie group
\[ G_0 = \left\{ \begin{bmatrix} A & b \\ 0 & c \end{bmatrix}: A \in
O(2),\ b \in \bb{R}^2,\ c \in \bb{R}^{\ast} \right\} . \]
We apply the method of equivalence to this $G_0$-structure, and after
two reductions we arrive at the bundle of {\em $2$-adapted}
coframings.  This is a $G_2$-structure $\scB_2 \to \scX$, where $G_2$
is the Lie group
\[ G_2 = \left\{ \begin{bmatrix} A & 0 \\ 0 & \det A \end{bmatrix}: A
\in O(2) \right\} . \]
There is a canonical coframing $(\w^1, \w^2, \w^3, \alp)$ (also known
as an {\em $(e)$-structure}) on $\scB_2$ whose structure equations
are
\begin{align}
d\w^1 & = -\alp \& \w^2 + A_1\, \w^2 \& \w^3 + A_2\, \w^3 \& \w^1 \notag \\
d\w^2 & = \alp \& \w^1 + A_2\, \w^2 \& \w^3 - A_1\, \w^3 \& \w^1
\label{subRiemstruct} \\
d\w^3 & = \w^1 \& \w^2 \notag \\
d\alp & = S_1\, \w^2 \& \w^3 + S_2\, \w^3 \& \w^1 + K\, \w^1 \& \w^2. \notag
\end{align}
Differentiating these equations shows that
\begin{align*}
dA_1 & = -2A_2\alp + \sum_{i=1}^3 B_{1i}\, \w^i \\
dA_2 & = 2A_1\alp + \sum_{i=1}^3 B_{2i}\, \w^i
\end{align*}
for some functions $B_{ij}$ on $\scB_2$, and that
\[ S_1  = B_{12} - B_{21}, \qquad S_2  = B_{11} + B_{22}. \]
By the general theory of $(e)$-structures, the functions $A_1, A_2,
K$ form a complete set of differential invariants for the
$G_2$-structure $\scB_2 \to \scX$, and hence for the sub-Riemannian
structure $(\scX, D, \langle, \rangle)$.

For later use, we observe that $\scB_2$ may be viewed geometrically
as a double cover of the unit circle bundle of the sub-Riemannian
metric.
If the sub-Riemannian structure $(\scX, D, \langle, \rangle)$ is {\em
orientable} (i.e., if we can choose an orientation on each of the
subspaces $D_x$ which varies smoothly on $\scX$), then $\scB_2$
consists of two disjoint connected components.  In this case
we can restrict the set of $0$-adapted coframings by requiring
that such a coframing be {\em oriented}, i.e., that the $2$-form
$\eta^1 \& \eta^2$ be a positive area form on $D$.  Doing so replaces
the $O(2)$ component of the structure group by $SO(2)$.  This does
not change anything essential in the preceding discussion, but it
does lead to a $G_2$-structure $\scB_2$ which is connected and is
naturally isomorphic to the unit circle bundle of $(\scX, D, \langle,
\rangle)$.

\section{Review of Finsler geometry of surfaces}

The material in this section is taken from \cite{Bryant96}.  (We
will, however, use the more standard notation for the invariants
which is found in \cite{BCS00}.)

A Finsler metric on a surface $\scM$ is determined by its indicatrix
bundle:  this is a smooth hypersurface $\Sigma^3 \subset T\scM$ with the
property that each fiber $\Sigma_x = \Sigma \cap T_x\scM$ is a smooth,
strongly convex curve which surrounds the origin $0_x \in T_x\scM$.  A
$3$-manifold $\Sigma \subset T\scM$ satisfying this condition is called
a {\em Finsler structure} on $\scM$.  A differentiable curve $\gamma:
[a,b] \to \scM$ is called a {\em $\Sigma$-curve} if, for every $s \in
[a, b]$, the velocity vector $\gamma'(s)$ lies in $\Sigma$.  The
following result is taken from \cite{Bryant96} and is due to Cartan
\cite{Cartan34}:

\begin{propvoid}
Let $\Sigma \subset T\scM$ be a Finsler structure on an oriented surface
$\scM$, with basepoint projection $\pi:\Sigma \to \scM$.  Then there exists
a unique coframing $(\w^1, \w^2, \alpha)$ on $\Sigma$ with the following
properties:
\begin{enumerate}
\item{$\w^1 \& \w^2$ is a positive multiple of any $\pi$-pullback of
a positive $2$-form on $\scM$.}
\item{The tangential lift $\gamma'$ of any $\Sigma$-curve satisfies
$(\gamma')^{\ast} \w^2 = 0$ and $(\gamma')^{\ast} \w^1 = dt$.}
\item{$d\w^1 \& \w^2 = d\w^2 \& \alpha = 0$.}
\item{$\w^1 \& d\w^1 = \w^2 \& d\w^2$.}
\item{$d\w^1 = -\alpha \& \w^2$.}
\end{enumerate}
Moreover, there exist functions $I, J, K$ on $\Sigma$ such that
\begin{align}
d\w^1 & = -\alpha \& \w^2 \notag \\
d\w^2 & = \alpha \& \w^1 - I\, \alp \& \w^2 \label{FSstructureeqs} \\
d\alpha & = K\, \w^1 \& \w^2 + J\, \alpha \& \w^2. \notag
\end{align}
\end{propvoid}

The Finsler structure on $\scM$ is Riemannian if and only if $I \equiv
0$; in this case, differentiating \eqref{FSstructureeqs} shows that
$J \equiv 0$ as well, and we recover the familiar structure equations
\begin{align*}
d\w^1 & = -\alpha \& \w^2 \\
d\w^2 & = \alpha \& \w^1 \\
d\alpha & = K\, \w^1 \& \w^2
\end{align*}
for an orthonormal coframing $(\w^1, \w^2)$ on $\scM$.  In this case, $\alp$ is
the Levi-Civita connection form, and $K$ is the usual Gauss curvature
on the surface.  For general Finsler surfaces, the function $K$
(called the {\em flag curvature}) is a well-defined function only on
$\Sigma$, not on $\scM$.

\section{The sub-Finsler equivalence problem}

Let $(\scX, D, F)$ be a sub-Finsler manifold consisting of a
three-dimensional manifold $\scX$, a rank two contact distribution 
$D$ on $\scX$,
and a sub-Finsler metric $F$ on $D$.  (Recall that $D$ is {\em
contact} if, for any two vector fields $\mbf{v}_1,
\mbf{v}_2$ locally spanning $D$, the vectors $\mbf{v}_1, \mbf{v}_2$, and
$[\mbf{v}_1, \mbf{v}_2]$ span the tangent space of $\scX$ at each point.)
As in the Finsler case, the sub-Finsler metric $F$ is completely
determined by its indicatrix bundle
\[ \Sigma = \{\mbf{u} \in D \mid F(\mbf{u}) = 1\}. \]
$\Sigma$ has dimension four, and each fiber $\Sigma_x = \Sigma \cap D_x$
is a smooth, strongly convex curve in $D_x$ which surrounds the
origin $0_x \in D_x$.  A 4-manifold $\Sigma \subset T\scX$ satisfying
this condition will be called a {\em sub-Finsler structure} on $(\scX,
D)$.

We will compute invariants for sub-Finsler structures via Cartan's
method of equivalence.  We begin by constructing a coframing on
$\Sigma$ which is nicely adapted to the sub-Finsler structure; this
procedure closely follows that used in \cite{Bryant95} for
constructing an adapted coframing for a Finsler structure on a
surface.

Let $g$ be any fixed sub-Riemannian metric on $(\scX, D)$, and let
$\Sigma_1$ be the unit circle bundle for $g$.  Then there exists a
well-defined, smooth function $r: \Sigma_1 \to \bb{R}^+$ with the
property that
\[ \Sigma = \{r(\mbf{u})^{-1}\, \mbf{u} \mid \mbf{u} \in \Sigma_1 \}. \]
Let $\rho: \Sigma \to \Sigma_1$ be the diffeomorphism which is the
inverse of the scaling map defined by $r$; i.e., $\rho$ satisfies
\[ \rho(r(\mbf{u})^{-1}\, \mbf{u}) = \mbf{u} \]
for $\mbf{u} \in \Sigma_1$.

Let $\pi:\Sigma \to \scX,\ \pi_1: \Sigma_1 \to \scX$ denote the respective
basepoint projections, and let $\mbf{u} \in \Sigma$.  (We trust that
using the same notation for points in $\Sigma$ and in $\Sigma_1$ will
not cause undue confusion.)  We will say that a vector $\mbf{v} \in
T_{\mbf{u}} \Sigma$ is {\em monic} if $\pi'(\mbf{u})(\mbf{v}) =
\mbf{u}$.  Since $\pi'(\mbf{u}):T_{\mbf{u}}\Sigma \to
T_{\pi(\mbf{u})} \scX$ is surjective with a one-dimensional kernel, the
set of monic vectors in $T_{\mbf{u}}\Sigma$ is an affine line.  A
nonvanishing 1-form $\theta$ on $\Sigma$ will be called {\em null} if
$\theta(\mbf{v}) = 0$ for all monic vectors $\mbf{v}$, and a 1-form
$\omega$ on $\Sigma$ will be called {\em monic} if $\omega(\mbf{v}) =
1$ for all monic vectors $\mbf{v}$.  The set of null 1-forms spans a
two-dimensional subspace of $T^{\ast}_{\mbf{u}}\Sigma$ at each point
$\mbf{u} \in \Sigma$, and the difference of any two monic 1-forms is
a null form.

In the sub-Riemannian case, $\omega^1$ is a monic form and the null
1-forms are spanned by $\omega^2$ and $\omega^3$.  (Recall that these
forms are part of the canonical coframing on $\Sigma_1$ described in
section \ref{subRiemsect}.)  Moreover, $D$ is defined by $D =
\{\omega^3\}^{\perp}$; this makes sense because according to the
structure equations \eqref{subRiemstruct}, $\omega^3$ descends to a
well-defined form on $\scX$.  Since the diagram
\setlength{\unitlength}{2pt}
\begin{center}
\begin{picture}(40,30)(0,0)
\put(5,25){\makebox(0,0){$\Sigma$}}
\put(35,25){\makebox(0,0){$\Sigma_1$}}
\put(20,5){\makebox(0,0){$\scX$}}
\put(8,20){\vector(3,-4){8}}
\put(31,20){\vector(-3,-4){8}}
\put(10,25){\vector(1,0){19}}
\put(7,15){\makebox(0,0){$\scriptstyle{\pi}$}}
\put(33,15){\makebox(0,0){$\scriptstyle{\pi_1}$}}
\put(20,28){\makebox(0,0){$\scriptstyle{\rho}$}}
\end{picture}
\end{center}
commutes, it is straightforward to verify that the null forms on
$\Sigma$ are spanned by $\rho^{\ast}(\omega^2)$ and
$\rho^{\ast}(\omega^3)$, that $D =
\{\rho^{\ast}(\omega^3)\}^{\perp}$, and that $\rho^{\ast}(r
\omega^1)$ is a monic form on $\Sigma$.

A local coframing $(\bar{\eta}^1, \bar{\eta}^2, \bar{\eta}^3,
\bar{\phi})$ on $\Sigma$ will be called {\em 0-adapted} if it
satisfies the conditions that $\bar{\eta}^1$ is a monic form,
$\bar{\eta}^2$ and $\bar{\eta}^3$ are null forms, and $D =
\{\bar{\eta}^3\}^{\perp}$.  For example, the coframing
\begin{equation}
\bar{\eta}^1 = \rho^{\ast}(r \omega^1), \qquad \bar{\eta}^2 =
\rho^{\ast}(\omega^2), \qquad \bar{\eta}^3 = \rho^{\ast}(\omega^3),
\qquad \bar{\phi} = \rho^{\ast}(\alp) \label{coframe0}
\end{equation}
is 0-adapted.  Any two 0-adapted coframings on $\Sigma$ vary by a
transformation of the form
\begin{equation}
\begin{bmatrix}
\tilde{\bar{\eta}}^1\\[0.1in] \tilde{\bar{\eta}}^2\\[0.1in]
\tilde{\bar{\eta}}^3\\[0.1in] \tilde{\bar{\phi}}
\end{bmatrix} =
\begin{bmatrix}
1 & a_1 & a_2 & 0 \\[0.1in] 0 & b_1 & b_2 & 0 \\[0.1in] 0 & 0 & b_3 &
0 \\[0.1in] c_1 & c_2 & c_3 & c_4
\end{bmatrix}^{-1}
\begin{bmatrix} \bar{\eta}^1\\[0.1in] \bar{\eta}^2\\[0.1in]
\bar{\eta}^3\\[0.1in] \bar{\phi}
\end{bmatrix}  \label{G0freedom}
\end{equation}
with $b_1 b_3 c_4 \neq 0$.  The set of all 0-adapted coframings
forms a principal fiber bundle $\scB_0 \to \Sigma$, with structure
group $G_0$ consisting of all matrices of the form \eqref{G0freedom}.
The right action of $G_0$ on sections $\sigma:\Sigma \to \scB_0$ is
given by $g\cdot
\sigma=g^{-1} \sigma$.  (This explains the inverse occurring in
\eqref{G0freedom}.)

There exist canonical 1-forms $\eta^1, \eta^2, \eta^3, \phi$ on
$\scB_0$ with the {\em reproducing property} that for any local
section $\sigma: \Sigma \to \scB_0$,
\[ \sigma^{\ast} (\eta^i) = \bar{\eta}^i, \qquad  \sigma^{\ast}
(\phi) = \bar{\phi}. \]
These are referred to as the {\em semi-basic} forms on $\scB_0$.
A standard argument shows that there also exist (non-unique) 1-forms $\alp_i,
\beta_i, \gamma_i$ (referred to as {\em pseudo-connection forms} or,
more succinctly, {\em connection forms}), linearly independent from
the semi-basic forms, and functions $T^i_{jk}$ on $\scB_0$ (referred
to as {\em torsion functions}) such that
\begin{equation}
\begin{bmatrix}
d\eta^1 \\[0.1in] d\eta^2 \\[0.1in] d\eta^3 \\[0.1in] d\phi
\end{bmatrix} =
-\begin{bmatrix}
0 & \alp_1 & \alp_2 & 0 \\[0.1in] 0 & \beta_1 & \beta_2 & 0 \\[0.1in]
0 & 0 & \beta_3 & 0 \\[0.1in] \gamma_1 & \gamma_2 & \gamma_3 &
\gamma_4
\end{bmatrix} \&
\begin{bmatrix}
\eta^1 \\[0.1in] \eta^2 \\[0.1in] \eta^3 \\[0.1in] \phi
\end{bmatrix} +
\begin{bmatrix}
T^1_{10}\, \eta^1 \& \phi \\[0.1in] T^2_{10}\, \eta^1 \& \phi
\\[0.1in] T^3_{12}\, \eta^1 \& \eta^2 \\[0.1in] 0
\end{bmatrix}. \label{G0structureeqs}
\end{equation}
These are the {\em structure equations} of the $G_0$-structure
$\scB_0$.  The semi-basic forms and connection forms together form a
local coframing on $\scB_0$.

We proceed with the method of equivalence by examining how the
functions $T^i_{jk}$ vary if we change from one 0-adapted coframing
to another.  A straightforward computation shows that under a
transformation of the form \eqref{G0freedom}, we have
\begin{align}
\tilde{T}^1_{10} & = c_4 T^1_{10} - \frac{a_1 c_4}{b_1} T^2_{10} \notag \\
\tilde{T}^2_{10} & = \frac{c_4}{b_1} T^2_{10} \label{G0torsiontrans} \\
\tilde{T}^3_{12} & = \frac{b_1}{b_3} T^3_{12}. \notag
\end{align}
In particular, the functions $T^2_{10}, T^3_{12}$ are {\em relative
invariants}: if they vanish for any 0-adapted coframing, then they
vanish for every 0-adapted coframing.  The coframing \eqref{coframe0}
has $T^2_{10}=-r^{-1},\ T^3_{12} = r^{-1}$, so we can assume that
these invariants are nonzero.  \eqref{G0torsiontrans} then implies
that we can adapt coframings to arrange that
\[ T^1_{10} = 0, \qquad T^2_{10}=-1, \qquad  T^3_{12} = 1. \]
A coframing satisfying this condition will be called {\em 1-adapted}.
For example, if we set
\[ dr = r_1\, \omega^1 + r_2\, \omega^2+ r_3\, \omega^3 + r_0\, \phi, \]
then the coframing
\begin{equation}
\bar{\eta}^1 = \rho^{\ast}(r \omega^1 - r_0\, \omega^2), \qquad
\bar{\eta}^2 = \rho^{\ast}(r \omega^2), \qquad \bar{\eta}^3 =
\rho^{\ast}(r^2 \omega^3), \qquad \bar{\phi} = \rho^{\ast}(\alp)
\label{coframe1}
\end{equation}
is 1-adapted.  Any two 1-adapted coframings on $\Sigma$ vary by a
transformation of the form
\begin{equation}
\begin{bmatrix}
\tilde{\bar{\eta}}^1\\[0.1in] \tilde{\bar{\eta}}^2\\[0.1in]
\tilde{\bar{\eta}}^3\\[0.1in] \tilde{\bar{\phi}}
\end{bmatrix} =
\begin{bmatrix}
1 & 0 & a_2 & 0 \\[0.1in] 0 & b_1 & b_2 & 0 \\[0.1in] 0 & 0 & b_1 & 0
\\[0.1in] c_1 & c_2 & c_3 & b_1
\end{bmatrix}^{-1}
\begin{bmatrix} \bar{\eta}^1\\[0.1in] \bar{\eta}^2\\[0.1in]
\bar{\eta}^3\\[0.1in] \bar{\phi}
\end{bmatrix}  \label{G1freedom}
\end{equation}
with $b_1 \neq 0$.  The set of all 1-adapted coframings forms a
principal fiber bundle $\scB_1 \subset \scB_0$, with structure
group $G_1$ consisting of all matrices of the form \eqref{G1freedom}.
When restricted to $\scB_1$, the connection forms $\alpha_1,\
\beta_3 - \beta_1, \ \gamma_4 - \beta_1$ become semi-basic, thereby
introducing new torsion terms into the structure equations of
$\scB_1$.  By adding multiples of the semi-basic forms to the 
connection forms so as to absorb as much of the torsion as possible, 
we can arrange that the structure equations of $\scB_1$ take the form
\begin{equation}
\begin{bmatrix}
d\eta^1 \\[0.1in] d\eta^2 \\[0.1in] d\eta^3 \\[0.1in] d\phi
\end{bmatrix} =
-\begin{bmatrix}
0 & 0 & \alp_2 & 0 \\[0.1in] 0 & \beta_1 & \beta_2 & 0 \\[0.1in] 0 &
0 & \beta_1 & 0 \\[0.1in] \gamma_1 & \gamma_2 & \gamma_3 & \beta_1
\end{bmatrix} \&
\begin{bmatrix}
\eta^1 \\[0.1in] \eta^2 \\[0.1in] \eta^3 \\[0.1in] \phi
\end{bmatrix} +
\begin{bmatrix}
T^1_{12}\, \eta^1 \& \eta^2 + T^1_{20}\, \eta^2  \& \phi \\[0.1in] -
\eta^1 \& \phi \\[0.1in] \eta^1 \& \eta^2 + T^3_{13}\, \eta^1 \&
\eta^3 + T^3_{30}\, \eta^3 \& \phi \\[0.1in] 0
\end{bmatrix}. \label{G1structureeqs}
\end{equation}
Moreover, we have
\begin{align*}
0 & \equiv d(d\eta^3) \mod{\eta^3} \\
& \equiv T^3_{30}\, \eta^1 \& \eta^2 \& \phi;
\end{align*}
therefore, $T^3_{30} = 0$.

We now repeat this process.  Under a transformation of the form
\eqref{G1freedom}, we have
\begin{align}
\tilde{T}^1_{20} & = b_1^2 T^1_{20} \notag \\
\tilde{T}^1_{12} & = b_1 T^1_{12} - b_1 c_1 T^1_{20} - a_2
\label{G1torsiontrans} \\
\tilde{T}^3_{13} & =  T^3_{13} + \frac{2b_2 + c_2}{b_1}. \notag
\end{align}
In particular, $T^1_{20}$ is a relative invariant which transforms by
a square, so its sign is fixed.  The coframing \eqref{coframe1} is
1-adapted, and if we set
\[ dr_0 = r_{01}\, \omega^1 + r_{02}\, \omega^2 + r_{03}\, \omega^3 +
r_{00}\, \phi, \]
it has $T^1_{20} = \dfrac{r + r_{00}}{r}$.  The condition that each
fiber of $\Sigma$ be a strongly convex curve enclosing the origin is exactly
the condition that this quantity be positive (see Lemma 
\ref{positivity} for a proof), so we can assume that
$T^1_{20} > 0$.  \eqref{G1torsiontrans} then implies that we can
adapt coframings to arrange that
\[ T^1_{20} = 1, \qquad T^1_{12} = T^3_{13} = 0. \]
A coframing satisfying this condition will be called {\em 2-adapted}.
Any two 2-adapted coframings on $\Sigma$ vary by a transformation of
the form
\begin{equation}
\begin{bmatrix}
\tilde{\bar{\eta}}^1\\[0.1in] \tilde{\bar{\eta}}^2\\[0.1in]
\tilde{\bar{\eta}}^3\\[0.1in] \tilde{\bar{\phi}}
\end{bmatrix} =
\begin{bmatrix}
1 & 0 & a_2 & 0 \\[0.1in] 0 & \eps & b_2 & 0 \\[0.1in] 0 & 0 & \eps &
0 \\[0.1in] -\eps a_2 & -2b_2 & c_3 & \eps
\end{bmatrix}^{-1}
\begin{bmatrix} \bar{\eta}^1\\[0.1in] \bar{\eta}^2\\[0.1in]
\bar{\eta}^3\\[0.1in] \bar{\phi}
\end{bmatrix}  \label{G2freedom}
\end{equation}
with $\eps = \pm 1$.  The set of all 2-adapted coframings forms a
principal fiber bundle $\scB_2 \subset \scB_1$, with structure
group $G_2$ consisting of all matrices of the form \eqref{G2freedom}.
When restricted to $\scB_2$, the connection forms $\beta_1,\
\gamma_1 + \alpha_2, \ \gamma_2 + 2\beta_2$ become semi-basic.  By 
adding multiples of the semi-basic forms to the connection forms so 
as to absorb as much of the torsion as possible, we can arrange that 
the structure equations of $\scB_2$ take the form
\begin{multline}
\begin{bmatrix}
d\eta^1 \\[0.1in] d\eta^2 \\[0.1in] d\eta^3 \\[0.1in] d\phi
\end{bmatrix} =
-\begin{bmatrix}
0 & 0 & \alp_2 & 0 \\[0.1in] 0 & 0 & \beta_2 & 0 \\[0.1in] 0 & 0 & 0
& 0 \\[0.1in] -\alp_2 & -2\beta_2 & \gamma_3 & 0
\end{bmatrix} \&
\begin{bmatrix}
\eta^1 \\[0.1in] \eta^2 \\[0.1in] \eta^3 \\[0.1in] \phi
\end{bmatrix} \\ +
\begin{bmatrix}
\eta^2 \& \phi \\[0.1in]
- \eta^1 \& \phi + T^2_{12}\, \eta^1 \& \eta^2 + T^2_{20}\, \eta^2 \&
\phi \\[0.1in]
\eta^1 \& \eta^2 + T^3_{23}\, \eta^2 \& \eta^3 - T^2_{12}\, \eta^3 \&
\eta^1 + T^2_{20}\, \eta^3 \& \phi \\[0.1in]
T^0_{12}\, \eta^1 \& \eta^2 + T^0_{10}\, \eta^1 \& \phi + T^0_{20}\,
\eta^2 \& \phi
\end{bmatrix}. \label{G2structureeqs}
\end{multline}
Moreover, we have
\begin{align*}
0 & \equiv d(d\eta^1) \mod{\eta^3} \\
& \equiv (T^2_{12} + T^0_{10})\, \eta^1 \& \eta^2 \& \phi;
\end{align*}
therefore, $T^0_{10} = -T^2_{12}$.

Under a transformation of the form \eqref{G2freedom}, we have
\begin{align}
\tilde{T}^2_{12} & = T^2_{12} + \eps(a_2 T^2_{20} + b_2)
\label{G2torsiontrans}\\
\tilde{T}^3_{23} & = \eps T^3_{23} + 2 b_2 T^2_{20} - a_2, \notag  \end{align}
so we can adapt coframings to arrange that
\[ T^2_{12} = T^3_{23} = 0. \]
A coframing satisfying this condition will be called {\em 3-adapted}.
Any two 3-adapted coframings on $\Sigma$ vary by a transformation of
the form
\begin{equation}
\begin{bmatrix}
\tilde{\bar{\eta}}^1\\[0.1in] \tilde{\bar{\eta}}^2\\[0.1in]
\tilde{\bar{\eta}}^3\\[0.1in] \tilde{\bar{\phi}}
\end{bmatrix} =
\begin{bmatrix}
1 & 0 & 0 & 0 \\[0.1in] 0 & \eps & 0 & 0 \\[0.1in] 0 & 0 & \eps & 0
\\[0.1in]  0 & 0 & c_3 & \eps
\end{bmatrix}^{-1}
\begin{bmatrix} \bar{\eta}^1\\[0.1in] \bar{\eta}^2\\[0.1in]
\bar{\eta}^3\\[0.1in] \bar{\phi}
\end{bmatrix} . \label{G3freedom}
\end{equation}
The set of all 3-adapted coframings forms a principal fiber bundle
$\scB_3 \subset \scB_2$, with structure group $G_3$ consisting of
all matrices of the form \eqref{G3freedom}.  When restricted to
$\scB_3$, the connection forms $\alpha_2,\ \beta_2$ become
semi-basic.  By adding multiples of the semi-basic forms to the 
connection forms so as to absorb as much of the torsion as possible, 
we can arrange that the structure equations of $\scB_3$ take
the form
\begin{multline}
\begin{bmatrix}
d\eta^1 \\[0.1in] d\eta^2 \\[0.1in] d\eta^3 \\[0.1in] d\phi
\end{bmatrix} =
-\begin{bmatrix}
0 & 0 & 0 & 0 \\[0.1in] 0 & 0 & 0 & 0 \\[0.1in] 0 & 0 & 0 & 0
\\[0.1in] 0 & 0 & \gamma_3 & 0
\end{bmatrix} \&
\begin{bmatrix}
\eta^1 \\[0.1in] \eta^2 \\[0.1in] \eta^3 \\[0.1in] \phi
\end{bmatrix} \\ +
\begin{bmatrix}
\eta^2 \& \phi + T^1_{13}\, \eta^1 \& \eta^3 + T^1_{23}\, \eta^2 \&
\eta^3 + T^1_{30}\, \eta^3 \& \phi \\[0.1in]
- \eta^1 \& \phi + T^2_{20}\, \eta^2 \& \phi + T^2_{13}\, \eta^1 \&
\eta^3 + T^2_{23}\, \eta^2 \& \eta^3 + T^2_{30}\, \eta^3 \&
\phi\\[0.1in]
\eta^1 \& \eta^2 + T^2_{20} \eta^3 \& \phi \\[0.1in]
T^0_{12}\, \eta^1 \& \eta^2 - T^1_{30}\, \eta^1 \& \phi + T^0_{20}\,
\eta^2 \& \phi
\end{bmatrix}. \label{G3structureeqs}
\end{multline}
(The coefficients $T^0_{12},\, T^0_{20}$ in \eqref{G3structureeqs}
are slightly modified from those in \eqref{G2structureeqs}.)
Moreover, we have
\begin{align*}
0 & \equiv d(d\eta^3) \mod{\phi} \\
& \equiv -(T^1_{13} + T^2_{23} + T^2_{20} T^0_{12})\, \eta^1 \&
\eta^2 \& \eta^3;
\end{align*}
therefore, we can set
\[ T^1_{13} = -\tfrac{1}{2}T^2_{20} T^0_{12} - A_2, \qquad T^2_{23} =
-\tfrac{1}{2}T^2_{20} T^0_{12} + A_2 \]
for some function $A_2$ on $\scB_3$.  (The reason for this choice of
notation will shortly become apparent.)

Under a transformation of the form \eqref{G3freedom}, we have
\begin{align}
\tilde{T}^1_{23} & = T^1_{23} + \eps c_3 \label{G3torsiontrans}\\
\tilde{T}^2_{13} & = T^2_{13} - \eps c_3 \notag  \end{align}
so we can adapt coframings to arrange that
\[ T^1_{23} = T^2_{13} = A_1 \]
for some function $A_1$.
A coframing satisfying this condition will be called {\em 4-adapted}.
Any two 4-adapted coframings on $\Sigma$ vary by a transformation of
the form
\begin{equation}
\begin{bmatrix}
\tilde{\bar{\eta}}^1\\[0.1in] \tilde{\bar{\eta}}^2\\[0.1in]
\tilde{\bar{\eta}}^3\\[0.1in] \tilde{\bar{\phi}}
\end{bmatrix} =
\begin{bmatrix}
1 & 0 & 0 & 0 \\[0.1in] 0 & \eps & 0 & 0 \\[0.1in] 0 & 0 & \eps & 0
\\[0.1in]  0 & 0 & 0 & \eps
\end{bmatrix}^{-1}
\begin{bmatrix} \bar{\eta}^1\\[0.1in] \bar{\eta}^2\\[0.1in]
\bar{\eta}^3\\[0.1in] \bar{\phi}
\end{bmatrix} . \label{G4freedom}
\end{equation}
The set of all 4-adapted coframings forms a principal fiber bundle
$\scB_4 \subset \scB_3$, with structure group $G_4 =
\bb{Z}/2\bb{Z}$.  $\scB_4$ is thus a double cover of $\Sigma$, and
the 1-forms $(\eta^1,\, \eta^2,\, \eta^3,\, \phi)$ form a canonical coframing
on $\scB_4$.  When restricted to $\scB_4$,
the last remaining connection form $\gamma_3$ becomes semi-basic, and
the structure equations of $\scB_4$ take the form
\begin{equation}
\begin{bmatrix}
d\eta^1 \\[0.1in] d\eta^2 \\[0.1in] d\eta^3 \\[0.1in] d\phi
\end{bmatrix} =
\begin{bmatrix}
\eta^2 \& \phi - (A_2 + \tfrac{1}{2} T^2_{20} T^0_{12})\, \eta^1 \&
\eta^3 + A_1\, \eta^2 \& \eta^3 + T^1_{30}\, \eta^3 \& \phi \\[0.1in]
- \eta^1 \& \phi + T^2_{20}\, \eta^2 \& \phi + A_1\, \eta^1 \& \eta^3
+ (A_2 - \tfrac{1}{2} T^2_{20} T^0_{12})\, \eta^2 \& \eta^3 +
T^2_{30}\, \eta^3 \& \phi\\[0.1in]
\eta^1 \& \eta^2 + T^2_{20} \eta^3 \& \phi \\[0.1in]
T^0_{12}\, \eta^1 \& \eta^2 + T^0_{13}\, \eta^1 \& \eta^3 +
T^0_{23}\, \eta^2 \& \eta^3 - T^1_{30}\, \eta^1 \& \phi + T^0_{20}\,
\eta^2 \& \phi + T^0_{30}\, \eta^3 \& \phi
\end{bmatrix}. \label{G4structureeqs}
\end{equation}
Finally, we differentiate equations \eqref{G4structureeqs} in order
to find any remaining relations among the torsion functions.  Setting
\[ dT^i_{jk} = T^i_{jk,1}\, \eta^1 + T^i_{jk,2}\, \eta^2 +
T^i_{jk,3}\, \eta^3 + T^i_{jk,0}\, \phi, \]
and computing $d(d\eta^3) = 0$ yields
\begin{align}
T^2_{20,1} & = T^2_{30} + T^2_{20} T^1_{30} \label{derivs} \\
T^2_{20,2} & = -(T^1_{30} + T^2_{20} T^0_{20}). \notag
\end{align}
Then computing $d(d\eta^2) \equiv 0 \mod{\eta^3}$ yields
\[ T^0_{20} = -2 T^2_{30}. \]
Further differentiation yields only differential equations for the
torsion functions and no new functional relations.

If we rename the $T^i_{jk}$ as follows:
\begin{align*}
T^2_{20} & = I \\
T^1_{30} & = J_1 \\
T^2_{30} & = J_2 \\
T^0_{12} & = K \\
T^0_{30} & = S_0 \\
T^0_{23} & = S_1 \\
T^0_{13} & = -S_2,
\end{align*}
then the structure equations on $\scB_4$ become
\begin{equation}
\begin{bmatrix}
d\eta^1 \\[0.1in] d\eta^2 \\[0.1in] d\eta^3 \\[0.1in] d\phi
\end{bmatrix} =
\begin{bmatrix}
\eta^2 \& \phi + A_1\, \eta^2 \& \eta^3 + (A_2 + \tfrac{1}{2} IK)\,
\eta^3 \& \eta^1 + J_1\, \eta^3 \& \phi \\[0.1in]
- \eta^1 \& \phi  + (A_2 - \tfrac{1}{2} IK)\, \eta^2 \& \eta^3 -
A_1\, \eta^3 \& \eta^1 + J_2\, \eta^3 \& \phi + I\, \eta^2 \&
\phi\\[0.1in]
\eta^1 \& \eta^2 + I\, \eta^3 \& \phi \\[0.1in]
S_0\, \eta^3 \& \phi + S_1\, \eta^2 \& \eta^3  + S_2\, \eta^3 \&
\eta^1  - J_1\, \eta^1 \& \phi - 2 J_2\, \eta^2 \& \phi + K\, \eta^1
\& \eta^2
\end{bmatrix}. \label{finalstructureeqs}
\end{equation}
(Compare with the sub-Riemannian structure equations
\eqref{subRiemstruct}.)

Our first result is that $I$ is the fundamental invariant that
determines whether or not a sub-Finsler structure is sub-Riemannian:

\begin{theorem}
The sub-Finsler structure $\Sigma$ is the unit circle bundle for a
sub-Riemannian metric if and only if $I \equiv 0$.
\end{theorem}
\begin{proof}
One direction is trivial: if $\Sigma = \Sigma_1$ for some
sub-Riemannian metric, then the canonical coframing which we have
constructed on $\Sigma$ is simply
\[ (\eta^1,\, \eta^2,\, \eta^3,\, \phi) = (\omega^1,\, \omega^2,\,
\omega^3,\, \alpha), \]
and so the structure equations \eqref{finalstructureeqs} must reduce
to \eqref{subRiemstruct}; therefore, $I \equiv 0$.

Now suppose that $I \equiv 0$.  Then
\[ 0 = d(d\eta^3) = (J_1\, \eta^2 - J_2\, \eta^1) \& \eta^3 \& \phi; \]
therefore, $J_1 \equiv J_2 \equiv 0$.  Now computing $d(d\eta^1) \equiv 0
\mod{\eta^1}$ shows that
\[ dA_1 \equiv (S_0 - 2 A_2)\, \phi \mod{\eta^1,\, \eta^2,\, \eta^3}, \]
and
computing $d(d\eta^2) \equiv 0 \mod{\eta^2}$ shows that
\[ dA_1 \equiv (-S_0 - 2 A_2)\, \phi \mod{\eta^1,\, \eta^2,\, \eta^3}. \]
Therefore, $S_0 \equiv 0$, and the structure equations
\eqref{finalstructureeqs} have the form \eqref{subRiemstruct}.  This
implies that $\Sigma$ is the unit circle bundle for a sub-Riemannian
metric, as desired.
\end{proof}

\section{The geodesic equations}

In this section we consider the problem of finding geodesics of the 
sub-Finsler structure.  Recall that the sub-Finsler length of a 
horizontal curve $\gamma:[a, b] \to \scX$ is given by
\begin{equation}\label{sf_length}
    \calL(\gamma) = \int_a^b F(\gamma'(t))\, dt.
\end{equation}
Finding critical points of this functional amounts to solving a 
constrained variational problem.  However, care must be taken when 
computing variations among horizontal curves on a non-integrable rank 
$s$ distribution $D$.  Given a horizontal curve $\gamma$, one would 
like to consider ``$D$-variational vector fields on $\gamma$ that 
vanish at the endpoints,'' but in general the existence of such 
vector fields is far from guaranteed.  In fact, this can fail 
spectacularly:  for example, when $D$ is an Engel system on a 
4-manifold $\scM$, $\scM$ is foliated by horizontal curves that have 
no such variations \cite{BH93}.

If such a vector field exists along $\gamma$, then $\gamma$ is said 
to be \emph{regular}, and the methods outlined in \cite{Griffiths83} 
suffice to find the first variation.  A horizontal curve for which 
this fails is called \emph{non-regular} (or \emph{abnormal}).  In 
\cite{Hsu92} Lucas Hsu established the following criterion for a 
curve to be non-regular:
\begin{theorem}(Hsu, \cite{Hsu92})
Let $\calI \subset T^*\scX$ be the annihilator of the rank $s$ 
distribution $D \subset T\scX$, and let $\Psi$ be the pullback of the 
canonical symplectic 2-form on $T^*\scX$ to $\calI$.  A horizontal 
curve $\gamma: [a,b] \rightarrow \scX$ is non-regular if and only if 
it has a lifting $\tilde{\gamma}: [a,b] \rightarrow \calI$ that does 
not intersect the zero section and satisfies $\tilde{\gamma}'(t)\hook 
\Psi = 0$ for all $t \in [a,b]$.
\end{theorem}

In the present case, $D$ is a contact system on a 3-manifold with 
$\calI = \mbox{span }\{\eta^3\}$, and it is easy to see that in this 
case all horizontal curves must be regular.  In what follows we will 
therefore use the variational methods described in 
\cite{Griffiths83}; our argument closely follows that of 
\cite{Hughen95}.

Choose an orientation of $D$, and consider the set of coframes in 
$\scB_{4}$ that preserve this orientation; for simplicity we will 
continue to use the notation $\scB_{4}$ for this set.  Every 
horizontal curve $\gamma:[a,b] \rightarrow \scX$ lifts to an integral 
curve $\bar{\gamma}: [a,b] \rightarrow \scB_{4}$ of the differential 
system $\bar{\calI} = \{ \eta^{2} , \eta^{3} \}$ with 
$\eta^1(\bar{\gamma}'(t)) \neq 0$.  This lift corresponds to choosing 
a 4-adapted coframing along the horizontal curve so that the vector 
$e_{1}$ dual to $\eta^1$ points in the direction of the velocity 
vector of the curve.   The sub-Finsler length of $\gamma$ is then 
equal to the integral of the monic one-form $\eta^1$ along the lifted 
curve $\bar{\gamma}$.   The problem of finding critical curves of the 
sub-Finsler length functional among horizontal curves is thus 
equivalent to finding critical curves of
\begin{equation}\label{sf_length2}
   \bar{\calL}(\bar{\gamma}) = \int_{\bar{\gamma}} \eta^1
\end{equation}
among integral curves $\bar{\gamma}$ of $\bar{\calI} = \{\eta^2, 
\eta^3 \}$ on $\scB_{4}$.

\begin{proposition}\label{P:first_var}  The critical curves of 
$\bar{\calL}$ among integral curves of $\bar{\calI}$ on $\scB_{4}$ 
are precisely the projections of integral curves, with transversality 
condition $\eta^1 \neq 0$, of the differential system $\calJ = \{ 
\eta^2, \eta^3, \phi - \lambda \eta^1, d\lambda - C\eta^1\}$ on $\scY 
\cong \scB_{4} \times \mathbb{R}$, where $\lambda$ is the coordinate 
on the fiber $\mathbb{R}$ and $C =  \lambda^{2} I + \lambda J_1+ A_2 
+ \frac{1}{2} IK$.
\end{proposition}

\begin{proof}
Following the algorithm in \cite{Griffiths83}, we define a 
submanifold $\scZ \subset T^*\scB_{4}$ as follows:  for each $x \in 
\scB_{4}$, let $\scZ_{x} = \eta^{1}(x) + \mbox{ span}\{ 
\bar{\calI}_x\}$
and let
\begin{equation}\label{Z_def}
   \scZ = \bigcup_{x \in \scB_{4}} \scZ_{x}.
\end{equation}
Let $\zeta$ be the pullback to $\scZ$ of the canonical 1-form on 
$T^*\scB_{4}$.  By the ``self-reproducing" property of $\zeta$, we 
may write
\begin{equation}\label{zeta_def}
   \zeta = \eta^1 + \lambda_{2} \eta^2 + \lambda_{3} \eta^3
\end{equation}
(where we have suppressed the obvious pullbacks in our notation). 
According to the general theory described in \cite{Griffiths83}, the 
critical points of the functional
\begin{equation}\label{zeta_functional}
   \mathcal{\tilde{L}}(\tilde{\gamma}) = \int_{\tilde{\gamma}} \zeta
\end{equation}
among unconstrained curves $\tilde{\gamma}$ on $\scZ$ project to 
critical curves of $\bar{\mathcal{L}}$ among integral curves 
$\bar{\gamma}$ of $\bar{\calI}$ on $\scB_{4}$; moreover, a curve 
$\tilde{\gamma}$ on $\scZ$ is a critical curve of 
$\mathcal{\tilde{L}}$ if and only if $\tilde{\gamma}'(t) \hook d\zeta 
|_{\tilde{\gamma}(t)} = 0$.

A straightforward computation shows that
\begin{align}
   d\zeta &= \lambda_{2}\, \phi \& \eta^1 - (1 + \lambda_{2} I)\phi \& 
\eta^2 - (\lambda_{3} I + J_{1} + \lambda_{2}J_{2})\phi \& \eta^3 + 
\lambda_{3}\eta^1\&\eta^2  \notag\\
    & \qquad + ( A_{2} + \tfrac{1}{2}IK - 
\lambda_{2}A_{1})\eta^3\&\eta^1 + (A_1 - \tfrac{1}{2} \lambda_2 IK + 
\lambda_2  A_2 ) \eta^2\&\eta^3 \label{dzeta} \\
& \qquad + d\lambda_{2} \& \eta^2 + d\lambda_{3} \& \eta^3.
\notag
\end{align}
By contracting $d\zeta$ with the vector fields dual to the coframing 
$\{\eta^1,\eta^2,\eta^3,\phi,d\lambda_{2},d\lambda_{3}\}$ on $\scZ$, 
we find that subject to the condition $\tilde{\gamma}^{*}\eta^1 \neq 
0$, the requirement that $\tilde{\gamma}' \hook d\zeta = 0$ is 
equivalent to the condition that $\tilde{\gamma}$ is an integral 
curve of the system
\[ \calJ = \{\eta^2, \eta^3, \phi - \lambda_{3} \eta^1, d\lambda_{3} 
- ( \lambda_{3}^{2} I + \lambda_{3} J_1 + A_2 + \tfrac{1}{2} IK) 
\eta^1 \}\]
on the submanifold $\scY \subset \scZ$ defined by $\lambda_{2} = 0$.
(Henceforth we will omit the subscript on $\lambda_{3}$.)  Curves 
satisfying this requirement project to critical curves of the 
functional $ \bar{\calL}$ among integral curves of $\bar{\calI}$ on 
$\scB_{4}$, and thus to local minimizers of the sub-Finsler length 
functional $\calL$ on $\scX$.  Since every horizontal curve on $\scX$ 
is regular, every local minimizer arises in this way.

\end{proof}

We will call a unit speed horizontal curve $\gamma : [a,b] \to \scM$ 
a \emph{sub-Finsler geodesic} if it has a lift to an integral curve 
of $\calJ$ on $\scY$.  When $\gamma$ has unit speed, it lifts to an 
integral curve of $\calJ$ if and only if it satisfies the 
\emph{geodesic equations}
\begin{equation}\label{geodesic_eqns}
    \eta^1 = ds,\, \eta^2 = 0,\, \eta^3 = 0,\, \phi = \lambda\, ds,\, 
d\lambda = C\, ds.
\end{equation}

\section{The Jacobi operator and the second variation}

This argument is similar to that given in \cite{Hughen95} for 
the sub-Riemannian case; we will describe it in some detail since 
\cite{Hughen95} is unpublished.

Since the geodesic equations are defined on the bundle $\scY \cong 
\scB_4 \times \bb{R}$, we will work on $\scY$ and use the
coframing $\{\eta^1, \eta^2, \eta^3, \eta^4, \eta^5\}$, where
\begin{align*}
\eta^4 & = \phi - \lambda \eta^1 \\
\eta^5 & = d\lambda - (\lambda^2 I+\lambda J_1 + A_2 + \tfrac{1}{2}IK)\,\eta^1.
\end{align*}
The structure equations \eqref{finalstructureeqs} imply that this 
coframing has structure equations
\begin{align}
d\eta^1 & = \eta^2 \& \eta^4 - \lambda\, \eta^1 \& \eta^2 + A_1\, 
\eta^2 \& \eta^3  + J_1\, \eta^3 \& \eta^4 -
(\lambda J_1 + A_2 + \tfrac{1}{2} IK)\, \eta^1 \& \eta^3 \notag \\
d\eta^2 & = - \eta^1 \& \eta^4 + J_2\, \eta^3 \& \eta^4   + I\, 
\eta^2 \& \eta^4 - \lambda I \, \eta^1 \& \eta^2 \notag
\\
& \qquad \qquad  + (A_2 - \tfrac{1}{2} IK)\, \eta^2 \& \eta^3 +
(A_1 - \lambda J_2 )\, \eta^1 \& \eta^3 \notag \\
d\eta^3 & = \eta^1 \& \eta^2 + I\, \eta^3 \& \eta^4 - \lambda I\, 
\eta^1 \& \eta^3 \label{Ystreqs} \\
d\eta^4 & = \eta^1 \& \eta^5 - J_1\, \eta^1 \& \eta^4 + (S_1 - 
\lambda A_1)\, \eta^2 \& \eta^3 + (S_0 - \lambda J_1)\,
\eta^3 \& \eta^4\notag \\
  & \qquad \qquad   - (\lambda + 2 J_2)\, \eta^2 \& \eta^4 + 
(\lambda^2  + 2\lambda J_2 + K)\, \eta^1 \& \eta^2
\notag
\\
& \qquad \qquad \qquad  + (J_1 \lambda^2 + (A_2 - S_0 + \tfrac{1}{2} 
IK)\lambda - S_2 ) \eta^1 \& \eta^3 \notag \\
d\eta^5 & =  - (\lambda^2 I+\lambda J_1 + A_2 + 
\tfrac{1}{2}IK)\,[\eta^2 \& \eta^4 - \lambda\, \eta^1 \& \eta^2 +
A_1\,\eta^2 \& \eta^3  + J_1\, \eta^3 \& \eta^4 \notag \\
& \qquad \qquad - (\lambda J_1 + A_2 + \tfrac{1}{2} IK)\, \eta^1 \& 
\eta^3] \notag \\
& \qquad \qquad \qquad+ \eta^1 \& d (\lambda^2 I+\lambda J_1 + A_2 +
\tfrac{1}{2}IK).
\notag
\end{align}

As in the previous section, every horizontal curve $\gamma$ has a 
canonical lift to an integral curve $\bar{\gamma}$ of the system 
$\bar{\calI} = \{\eta^2,\, \eta^3\}$ on $\scB_4$.  This in turn has a 
canonical lift to an integral curve $\tilde{\gamma}$ of the system 
$\tilde{\calI} = \{\eta^2,\, \eta^3,\, \eta^4\}$ on $\scY$.  The 
length of $\gamma$ is equal to the integral of $\eta^1$ along the 
lifted curve $\tilde{\gamma}$, and $\gamma$ is a geodesic if and only 
if $\tilde{\gamma}$ is an integral curve of the system $\calJ = 
\{\eta^2,\, \eta^3,\, \eta^4\,\, \eta^5\}$ on $\scY$.

Suppose that $\gamma:[0, \ell] \to \scX$ is a horizontal curve 
joining points $p$ and $q$ in $\scX$.  If $\gamma_t$ is a 
fixed-endpoint variation of $\gamma$ through horizontal curves, then 
$\gamma_t$ lifts to a variation $\tilde{\gamma}_t$ of 
$\tilde{\gamma}$ through integral curves of $\tilde{\calI}$; this 
variation does not necessarily fix endpoints, but it satisfies the 
condition
\[ \pi \circ \tilde{\gamma}_t(0) = p, \qquad \pi \circ 
\tilde{\gamma}_t(\ell) = q, \]
where $\pi: \scY \to \scX$ is the usual base point projection.  A 
variation $\tilde{\gamma}_t$ satisfying these conditions will be 
called an {\em admissible variation} of $\tilde{\gamma}$, and its 
variational vector field $\frac{\partial \tilde{\gamma}_t}{\partial 
t}$ at $t=0$ will be called an {\em infinitesimal admissible 
variation} along $\tilde{\gamma}$.

Now suppose that $\gamma$ is a geodesic.  Let $\gamma_{t,u}$ be 
2-parameter fixed-endpoint variation of $\gamma$ through horizontal 
curves, and let $\tilde{\gamma}_{t,u}$ be its lift to $\scY$, with 
infinitesimal admissible variations
\[ V(s) = \frac{\partial}{\partial t} \tilde{\gamma}_{t,u}(s) 
\mid_{t=u=0}, \qquad
W(s) = \frac{\partial}{\partial u} \tilde{\gamma}_{t,u}(s) \mid_{t=u=0}. \]
Let $(e_1, \ldots, e_5)$ be the framing dual to the coframing 
$(\eta^1, \ldots, \eta^5)$ on $\scY$, and write
\[ V(s) = \sum_{i=1}^5 V_i(s) e_i(s), \qquad W(s) = \sum_{i=1}^5 
W_i(s) e_i(s). \]
The Hessian $\calL_{**}(V,W)$ of the length functional is, by definition,
\[ \calL_{**}(V,W) = \frac{\partial^2}{\partial t \partial u} 
\calL(\gamma_{t,u})\mid_{t=u=0}. \]

\begin{proposition}
Let $\tilde{\gamma}:[0, \ell] \to \scY$ be a lifted geodesic with 
2-parameter admissible variation $\tilde{\gamma}_{t,u}$.  The Hessian 
of the length functional evaluated at the infinitesimal admissible 
variations $V = \sum V_i e_i$ and $W = \sum W_i e_i$ is
\[ \calL_{**}(V,W) = \int_0^{\ell} W_3\, J(V_3)\, ds, \]
where $J$ is a self-adjoint, fourth-order differential operator on 
the space of smooth functions on $[0, \ell]$ given by
\begin{equation}
J(u) = \ddddot{u} + \frac{d}{ds}\left(P \dot{u}\right) + Qu, \label{Jacobi}
\end{equation}
for certain functions $P,Q$ along $\tilde{\gamma}$.
\end{proposition}
Here the dots over $u$ represent derivatives with respect to $s$, and 
the precise definitions of $P$ and $Q$ will appear in the proof.

\begin{pf}
For any smooth function $f$ on $\scY$, we will write
\[ df = \sum_{i=1}^5 f_{,i} \eta^i. \]
Differentiating \eqref{Ystreqs} yields relations among the 
derivatives of the invariants of the sub-Finsler structure, and these 
relations must be taken into account in the computations that follow.

As in \cite{Hughen95}, the Hessian $\calL_{**}(V,W)$ is equal to the integral
\[ \int_{\tilde{\gamma}} W \hook d(V \hook d\zeta), \]
where $\zeta = \eta^1 + \lambda\, \eta^3$.  A long but 
straightforward computation shows that along $\tilde{\gamma}$, the 
integrand $W \hook d(V \hook d\zeta)$ is equivalent modulo $\calJ$ to
\begin{equation}
\begin{split}
\Big[ & W_2 \Big( \dot{V}_{4} + A_1 \dot{V}_{3} - (\lambda^2 + 2 J_2 
\lambda + A_1 + K) V_2 \\
& \qquad \qquad + \big(-J_1 \lambda^2 + (A_1 I - \tfrac{1}{2}IK - A_2 
+ S_0)\lambda + S_2\big)V_3 + J_1 V_4 - V_5 \Big) \\
& + W_3 \Big(-\dot{V}_{5} - A_1 \dot{V}_{2} + (I\lambda + J_1)\dot{V}_{4} \\
& \qquad \qquad + \big(-J_1 \lambda^2 + (A_1 I - \tfrac{1}{2}IK - A_2 
+ S_0)\lambda + S_2 - A_{1,1}\big)V_2 \\
& \qquad \qquad + \big((S_0 I + I_{,3} + S_{0,4})\lambda^2 + ((A_1 + 
K)J_2 + J_1 S_0 - S_1 - S_{0,1} + S_{2,4})\lambda \\
& \qquad \qquad \qquad + (\tfrac{1}{4}I^2K^2 + \tfrac{1}{2}(I K_{,3} 
+ I_{,3}K) + A_1^2 + A_2^2 + A_2 IK + J_1 S_2 + A_{2,3})\big) V_3 \\
& \qquad \qquad + \big(I_{,4} \lambda^2 + (IJ_1 + J_{1,4})\lambda + 
\tfrac{1}{2} (I^2 K + K I_{,4} - I A_{1,4})\\
& \qquad \qquad \qquad  - A_1 I^2 - A_2 I - 2 S_0 I + J_1^2 - A_1 + 
A_{2,4}\big) V_4 \\
& \qquad \qquad + I\lambda V_5 \Big) \\
& + W_4 \Big(-\dot{V}_{2} - (I\lambda + J_1)\dot{V}_{3} + J_1 V_2 + 
(-I^2\lambda^2 - (I J_1 + J_2)\lambda + A_1) V_3 - V_4 \Big) \\
& + W_5 \Big(\dot{V}_{3} - V_2 + I\lambda V_3 \Big)
\Big] \eta^1 .
\end{split} \label{integrand}
\end{equation}
(This computation takes into account the fact that along 
$\tilde{\gamma}, \ V_{i,1} = \dot{V}_i$.)

Now define $\Gamma(t,u,s) = \tilde{\gamma}_{t,u}(s)$.  Since each 
curve $\tilde{\gamma}_{t,u}$ is an integral curve of $\tilde{\calI}$, 
we have
\[ \Gamma^* \left[ \begin{array}{l} \eta^1 \\[0.05in] \eta^2 
\\[0.05in] \eta^3 \\[0.05in] \eta^4 \\[0.05in] \eta^5  \end{array} 
\right] =
\left[ \begin{array}{l}
V_1(t,u,s)\, dt + W_1(t,u,s)\, du + Y_1(t,u,s)\, ds \\[0.05in]
V_2(t,u,s)\, dt + W_2(t,u,s)\, du  \\[0.05in]
V_3(t,u,s)\, dt + W_3(t,u,s)\, du  \\[0.05in]
V_4(t,u,s)\, dt + W_4(t,u,s)\, du  \\[0.05in]
V_5(t,u,s)\, dt + W_5(t,u,s)\, du + Y_5(t,u,s)\, ds  \end{array} \right]  \]
for some functions $V_i, W_i, Y_i$ satisfying (with a minor abuse of 
notation) $V_i(0,0,s) = V_i(s)$, $W_i(0,0,s) = W_i(s)$.  Since 
$\tilde{\gamma}$ is the lift of a geodesic, we also have
\[ Y_1(0,0,s) = 1, \qquad Y_5(0,0,s) = 0, \qquad \frac{\partial 
Y_1}{\partial t}(0,0,s) = \frac{\partial Y_1}{\partial u}(0,0,s) = 0. 
\]
When the structure equations \eqref{Ystreqs} are pulled back by 
$\Gamma$ and then restricted to $\tilde{\gamma}$, they imply that
\begin{align}
\dot{V}_1 & = -\lambda V_2 - (J_1\lambda + A_2 + \tfrac{1}{2} IK) V_3 \notag \\
\dot{V}_2 & = -I\lambda V_2 + (-J_2\lambda + A_1) V_3 - V_4 \label{Vdiffeqs} \\
\dot{V}_3 & = V_2 - I\lambda V_3 \notag \\
\dot{V}_4 & = (\lambda^2 + 2J_2\lambda + K) V_2 + (J_1\lambda^2 + 
(A_2 + \tfrac{1}{2} IK - S_0)\lambda - S_2) V_3 - J_1 V_4 + V_5. 
\notag
\end{align}
The third equation in \eqref{Vdiffeqs} implies that $V_2 = \dot{V_3} 
+ I\lambda V_3$.  The first equation in \eqref{Vdiffeqs} can then be 
written as
\[ \frac{d}{ds} (V_1 + \lambda V_3) = 0. \]
Since $V$ is an admissible infinitesimal variation, $V_1, V_2$, and 
$V_3$ must vanish at the endpoints of $\tilde{\gamma}$, and it 
follows that $V_1 = -\lambda V_3$.  Equations \eqref{Vdiffeqs} can 
now be used to express $V_1, V_2, V_4$, and $V_5$ in terms of $V_3$ and its derivatives, as follows:
\begin{align}
V_1 & = -\lambda V_3 \notag \\
V_2 & = \dot{V}_3 + I\lambda V_3 \notag \\
V_4 & = -\ddot{V}_3 - 2I\lambda \dot{V}_3 + [-(2 I^2 + I_{,4}) 
\lambda^2 - (2 I J_1 + 2 J_2)\lambda + (A_1 - I A_2  - \tfrac{1}{2} 
I^2 K )] V_3 \notag \\
V_5 & =  -\dddot{V}_3 - (2 I \lambda + J_1) \ddot{V}_3 \label{Veqs} \\
& \qquad  + [-(4 I^2  + 3 I_{,4} + 1) \lambda^2 - (8 I J_1 + 6 J_2) 
\lambda + (A_1 - 3 I A_2 - \tfrac{3}{2} I^2 K - K )] \dot{V}_3 \notag 
\\
& \qquad + [-(4 I^3 + 6 I I_{,4} + I)\lambda^3 \notag \\
& \qquad \ \ \ \ \ \   - (12 I^2 J_1 + 5 I_{,4} J_1 + I_{,41} + 2 I 
J_{1,4}  + 8 I J_2 + 2 J_{2,4} + J_1)\lambda^2 \notag \\
& \qquad \ \ \ \ \ \   + (-3 I K I_{,4} - 3 A_2 I_{,4} + 2 A_{1,4} + 
I^2 A_{1,4} - 2 I A_{2,4} + 4 A_1 I  + 2 I^3 A_1 \notag \\
& \qquad \qquad \qquad   + 3 A_2 - 2 I^3 K  - 8 I J_1^2 - 8 J_1 J_2 - 
\tfrac{3}{2} I K + 3 S_0 + 4 I^2 S_0) \lambda  \notag \\
& \qquad \ \ \ \ \ \  + (A_{1,1} - I A_{2,1}- \tfrac{1}{2} I^2 K_{,1} 
+ A_1 J_1 - 4I A_2 J_1 - 3 A_2 J_2 - \tfrac{5}{2} I^2 J_1 K \notag \\
& \qquad \qquad \qquad  - 2 I J_2 K  + S_2)  ]  V_3.  \notag
\end{align}
When these expressions are substituted into the integrand 
\eqref{integrand}, it takes the form
\[ W_3 (\ddddot{V}_3 + \frac{d}{ds}(P \dot{V}_3) + Q V_3), \]
where
\begin{align*}
\scriptstyle P &
\scriptstyle = (2I^2 + 4 I_{,4} + 1)\lambda^2  + 8(J_1 I + 
J_2)\lambda + \left( 2 K I^2 + 4 A_2 I + K - 3 A_1 \right)   \\
\scriptstyle Q &
\scriptstyle = [4 I^4 + 2(8 I_{,4} + 1)I^2 + (5 I_{,4}^2 + I_{,4})] 
\lambda^4 \\
& \scriptstyle \ \  + [20 J_1 I^3 + (14 J_2 + 10 J_{1,4}) I^2 + (40 
J_1 I_{,4} + 6 I_{,41} + 4 J_1 + 8 J_{2,4}) I + (12 J_2 I_{,4} + 6 
J_{1,4} I_{,4} + J_2 + J_{1,4}) ] \lambda^3 \\
& \scriptstyle \ \ + [(3 K - 8 A_1) I^4 - (10 A_2 + 16 S_0 + 4 
A_{1,4}) I^3 + (16 K I_{,4} - 10 A_1 I_{,4} - 19 A_1 + 8 A_{2,4} + 40 
J_1^2 + K ) I^2 \\
& \scriptstyle \ \ \ \qquad + (5 A_2 I_{,4} - 4  A_{1,4} I_{,4} - 16 
S_0 I_{,4} + 42 J_1 J_2 + 18 J_1 J_{1,4} + 2 J_{1,41} - 12 A_2 - 
\tfrac{13}{2} A_{1,4} - 6 S_0 ) I \\
& \scriptstyle \ \ \ \qquad + (5 K I_{,4}^2 + 18 J_1^2 I_{,4} - 13 
A_1 I_{,4} + 6 A_{2,4} I_{,4} + 2 K I_{,4} + 6 J_1 I_{,41} + I_{,411} 
+ I_{,3} - 2 A_1  + A_{2,4} \\
& \scriptstyle \ \ \ \ \qquad \qquad  + 3 J_1^2 + 8 J_2^2 + 8 J_2 
J_{1,4} + 10 J_1 J_{2,4} + 2 J_{2,41} + S_{0,4})] \lambda^2 \\
& \scriptstyle \ \ + [(12 K J_1 - 20 A_1 J_1 - 2 A_{1,1} + K_{,1})I^3 \\
& \scriptstyle \ \ \ \qquad + (4 K J_{1,4} - 9 J_1 A_{1,4} - 2 
A_{2,1} - A_{1,41} - 14 J_1 A_2 - 12 J_2 A_1 + 9 J_2 K - 36 J_1 S_0 - 
4 S_{0,1}) I^2 \\
& \scriptstyle \ \ \ \qquad + (20 K J_1 I_{,4} + 4 K_{,1} I_{,4} + 4 
K I_{,41} + 7 A_2 J_{1,4} + 4 K J_{2,4} + 16 J_1 A_{2,4} - 5 J_2 
A_{1,4} - 3 A_{1,1} + 2 A_{2,41} \\
& \scriptstyle \ \ \ \ \qquad \qquad + \tfrac{3}{2} K_{,1} - 34 A_1 
J_1 - 4 A_2 J_2 + 4 K J_1 + 24 J_1^3 - 20 J_2 S_0 + 2 S_2) I \\
& \scriptstyle \ \ \ \qquad + (14 J_1 A_2 I_{,4} + 4 A_{2,1} I_{,4} + 
8 K J_2 I_{,4} + 5 A_2 I_{,41} + 7 A_2 J_{2,4} - 8 J_1 A_{1,4} + 8 
J_2 A_{2,4} - 3 A_{2,1} - 2 A_{1,41} \\
& \scriptstyle \ \ \ \ \qquad \qquad - 4 S_{0,1} + S_{2,4} - 18 A_1 
J_2- 13 A_2 J_1 + 24 J_1^2 J_2 + \tfrac{5}{2} K J_2 - 8 J_1 S_0 - 
S_1)] \lambda \\
& \scriptstyle \ \ + [(\tfrac{3}{4} K^2 - 3 A_1 K) I^4 - 
(\tfrac{3}{2} K A_{1,4} + 5 A_1 A_2 + 3 A_2 K + 6 K S_0) I^3 \\
& \scriptstyle \ \ \ \qquad + (\tfrac{5}{2} K^2 I_{,4} - \tfrac{5}{2} 
A_2 A_{1,4}  + 3 K A_{2,4} + 3 J_1 K_{,1} + \tfrac{1}{2} K_{,11} - 7 
A_2^2 + K^2 - \tfrac{13}{2} A_1 K + 9 K J_1^2 - 10 A_2 S_0 ) I^2 \\
& \scriptstyle \ \ \ \qquad + (6 A_2 K I_{,4} + 5 A_2 A_{2,4} + 4 J_1 
A_{2,1} - 3 K A_{1,4} + A_{2,11}+ 3 J_2 K_{,1} + \tfrac{1}{2} K_{,3} 
- 11 A_1 A_2 + 12 A_2 J_1^2 \\
& \scriptstyle \ \ \ \ \qquad \qquad - 3 A_2 K + 11 J_1 J_2 K- 
\tfrac{7}{2} K S_0) I \\
& \scriptstyle \ \ \ \qquad + (3 A_2^2 I_{,4}+ \tfrac{1}{2} K I_{,3} 
- 5 A_2 A_{1,4} + 4 J_2 A_{2,1} + A_{2,3} - A_{1,11} - S_{2,1} + 2 
A_1^2 - 8 A_2^2 + 12 A_2 J_1 J_2 \\
& \scriptstyle \ \ \ \ \qquad \qquad - 6 A_2 S_0 + 2 J_2^2 K + J_1 S_2)]
      .
\end{align*}
Thus the proposition is proved.
\end{pf}

We saw in the proof of this proposition that any infinitesimal 
admissible variation $V = \Sigma_i V_i e_i$ satisfies
\[ \dot{V}_3 = V_2 -  I \lambda V_3, \]
and that $V_1, V_2, V_3$ vanish at the endpoints 
$\tilde{\gamma}(0),\, \tilde{\gamma}(\ell)$ of $\tilde{\gamma}$.  In 
particular, $V_3$ vanishes to first order at $0$ and $\ell$.  Let 
$C^{\infty}_0[0, \ell]$ denote the space of smooth functions on 
$[0,\ell]$ that vanish to first order at the endpoints, and note that 
the Jacobi operator $J$ is formally self-adjoint on $C^{\infty}_0[0, 
\ell]$.

Define a quadratic form $\calQ(u)$ by
\[ \calQ(u) = \calL_{**}(u,u). \]
Recall that the {\em index} of $\calQ$ is the dimension of the 
largest subspace of $C^{\infty}_0[0, \ell]$ on which $\calQ$ is 
negative definite.  Because $J$ is self-adjoint on $C^{\infty}_0[0, 
\ell]$, its eigenvalues form a countable subset of the real numbers 
with $+\infty$ as the only possible cluster point.  It follows that 
$J$ has only finitely many negative eigenvalues, and that therefore 
the index of $\calQ$ is finite.

\begin{definition}
A point $c \in (0,\ell)$ is a {\em conjugate point of $J$ with 
multiplicity $m$} if the space of solutions of the system
\begin{equation}
  J(u) = 0, \qquad u(0) = \dot{u}(0) = 0 \label{JacobianIVP}
\end{equation}
which vanish to first order at $c$ has dimension $m > 0$.  The point 
$\gamma(c)$ along a geodesic $\gamma$ is a {\em conjugate point} of 
$\gamma$ if $c$ is a conjugate point of the corresponding Jacobi 
operator $J$ along $\tilde{\gamma}$.
\end{definition}
Note that, since $J$ is a fourth-order operator, the multiplicity of 
any conjugate point of $J$ is either one or two.

\begin{theorem}\label{JacobiThm}
The index of $\calQ$ is equal to the number of conjugate points of 
$J$, counted with multiplicity.
\end{theorem}

\begin{pf}
Suppose that the index of $\calQ$ is $n$, and let
\[ \lambda_1 \leq \lambda_2 \leq \cdots \leq \lambda_n < 0 \]
be the negative eigenvalues of $J$.  For each $s \in (0, \ell]$, let
\[ \Lambda_1(s) \leq \Lambda_2(s) \leq \cdots \]
denote the eigenvalues of the operator $J$ on the space $[0, s]$. 
(Note that $\Lambda_i(\ell) = \lambda_i$ for $i = 1, \ldots, n$.)  It 
follows from general theory (see, e.g., \cite{Kato66}) that each 
$\Lambda_i(s)$ is a strictly decreasing, continuous function on $(0, 
\ell]$, with $\lim_{s \to 0^+} \Lambda_i(s) = +\infty$.  Therefore, 
each function $\Lambda_i(s), \ i=1, \ldots, n$ has exactly one root 
$c_i$.

These roots
\[ 0 < c_1 \leq c_2 \leq \cdots \leq c_n < \ell \]
are precisely the conjugate points of $J$ between $0$ and $\ell$.  To 
see this, note that by definition, the condition $\Lambda_i(c_i) = 0$ 
implies that there exists a function $u_i \in C^{\infty}_0[0, c_i]$ 
satisfying $J(u_i) = 0$.  By extending $u_i$ to a solution of $J(u) = 
0$ on $C^{\infty}_0[0, \ell]$, we get a solution to 
\eqref{JacobianIVP} which vanishes to first order at $c_i$. 
Conversely, if $c \in (0, \ell)$ is a conjugate point, then $J$ has a 
zero eigenvalue on $C^{\infty}_0[0, c_i]$; therefore, $\Lambda_j(c) = 
0$ for some positive integer $j$.  Since $\Lambda_j$ is a strictly 
decreasing function, it follows that $\Lambda_j(\ell)$ is equal to 
one of the negative eigenvalues $\lambda_j$ of $J$, and therefore, $c 
= c_j$.
\end{pf}

\begin{corollary}
A geodesic no longer minimizes length beyond its first conjugate point.
\end{corollary}

\begin{pf}
Suppose that $\gamma(c)$ is the first conjugate point of the geodesic 
$\gamma$.  Theorem \ref{JacobiThm} implies that for every $\ell > c$, 
the index of $\calQ$ on the space $C^{\infty}_0[0, \ell]$ is 
positive, and so there exists a function $u \in C^{\infty}_0[0, 
\ell]$ for which $\calQ(u) < 0$.  Setting $V_3 = u$ and defining 
$V_1, V_2, V_4, V_5$ according to equations \eqref{Veqs} defines a 
direction $V = \sum_{i=1}^5 V_i e_i$ along which the length 
functional $\tilde{\calL}$ (and hence $\calL$) decreases.
\end{pf}

\section{Symmetries and homogeneous examples}

In this section we examine the symmetries of sub-Finsler structures
and describe homogeneous examples.

\begin{definition}
Let $\Sigma$ be a sub-Finsler structure on $(\scX,D)$.  A {\em symmetry}
of $\Sigma$ is a diffeomorphism $\Phi:\scX \to \scX$ which satisfies
$\Phi'(\Sigma) = \Sigma$.  (Note that this implies that $\Phi'(D) =
D$ as well.)  A {\em symmetry} of the $\bb{Z}/2\bb{Z}$-structure
$\scB_4$ is a diffeomorphism $\Psi:\scB_4 \to \scB_4$ with the
property that
\[ (\Psi^{\ast} \eta^1,\, \Psi^{\ast}\eta^2,\, \Psi^{\ast}\eta^3,\,
\Psi^{\ast}\phi) = (\eta^1,\, \eta^2,\, \eta^3,\, \phi). \]
\end{definition}
A standard argument shows that the map
\[ \Phi \to \Phi' \]
gives a one-to-one correspondence between orientation-preserving
symmetries of $\Sigma$ and symmetries of $\scB_4$.  By a theorem of
Kobayashi \cite{Kobayashi57}, it follows that the group of symmetries
of $\Sigma$ can be given the structure of a Lie group of dimension at
most four.

There are two possible definitions of homogeneity for a sub-Finsler
structure: we could say that $\Sigma \subset T\scX$ is homogeneous if
its group of symmetries acts transitively on $\scX$, or we could require
the more restrictive condition that this group act transitively on
$\Sigma$.  Both notions are interesting, and we will consider each of
them in the remainder of this section.

\subsection{Symmetry groups of dimension four}

First we consider the case where the group of symmetries of $\Sigma$
is four-dimensional and acts transitively on $\scB_4$.  Since any
symmetry must preserve the canonical coframing $(\eta^1,\, \eta^2,\,
\eta^3,\, \phi)$ on $\scB_4$, it follows that all the torsion
functions must be constants. Conversely, if all the torsion
coefficients are constants, then the structure equations of $\scB_4$
define a local Lie group structure on $\scB_4$ for which the
canonical coframing is left-invariant; this Lie group then acts
transitively on $\scB_4$ in the obvious way.

So, suppose that all the torsion functions in the structure equations
\eqref{finalstructureeqs} are constant.  Then
\[ 0 = d(d\eta^3) = [-(I J_1 + J_2)\, \eta^1 + (J_1 - 2 I J_2)\,
\eta^2] \& \eta^3 \& \phi; \]
therefore, $J_1 = J_2 = 0$.  Next we have
\begin{align*}
0 = d(d\eta^1) & \equiv -S_2\, \eta^1 \& \eta^2 \& \eta^3 \mod{\phi} \\
0 = d(d\eta^2) & \equiv S_1\, \eta^1 \& \eta^2 \& \eta^3 \mod{\phi};
\end{align*}
therefore, $S_1 = S_2 = 0$.  Then
\[ 0 = d(d\phi) \equiv (S_0 - IK)\, \eta^1 \& \eta^2 \& \phi \mod{\eta^3}; \]
therefore, $S_0 = IK$.  Now
\[ 0 = d(d\eta^1) = [(2A_1 + I A_2 + \tfrac{1}{2} I^2 K)\, \eta^1 +
(-2I A_1 + 2 A_2 - IK)\, \eta^2] \& \eta^3 \& \phi; \]
therefore,
\[ A_1 = -\frac{I^2 K}{I^2 + 2}, \qquad A_2  = -\frac{IK (I^2 - 2)}{2
(I^2 + 2)}. \]
Finally, we have
\[ 0 = d(d\eta^2) = \frac{4 IK}{I^2 + 2}\, \eta^1 \& \eta^3 \& \phi; \]
therefore, $IK = 0$.  If $I=0$, then $\Sigma$ is sub-Riemannian; the
homogeneous sub-Riemannian structures are classified in
\cite{Hughen95}.  So suppose that $I \neq 0$.  Then we have $K=0$, and the
structure equations \eqref{finalstructureeqs} reduce to
\begin{align}
d\eta^1 & = \eta^2 \& \phi \notag \\
d\eta^2 & = -\eta^1 \& \phi + I\, \eta^2 \& \phi \label{homogstruct} \\
d\eta^3 & = \eta^1 \& \eta^2 + I\, \eta^3 \& \phi \notag \\
d\phi & = 0. \notag
\end{align}
Differentiating \eqref{homogstruct} yields no additional 
restrictions;
thus there exists (at least locally) a 1-parameter family of 
homogeneous sub-Finsler structures which are not sub-Riemannian.

In fact, these structure equations can be integrated explicitly. 
First, since $d\phi=0$, there exists a function $\theta$ on $\Sigma$ 
such that
\[ \phi = d\theta. \]
Next, note that the system $S = \{\eta^1,\, \eta^2\}$ is Frobenius 
(i.e., $dS \equiv 0 \mod{S}$); therefore, there exist functions $x,y$ 
on $\Sigma$, independent from $\theta$, such that
\[ S = \{dx,\, dy\}, \]
and functions $a_{ij}, \ i,j=1,2,$ such that
\begin{align*}
\eta^1 & = a_{11}\, dx + a_{12}\, dy \\
\eta^2 & = a_{21}\, dx + a_{22}\, dy.
\end{align*}
Now the first two equations of \eqref{homogstruct} imply that the 
$a_{ij}$ are functions of $x,y,\theta$ alone, and that
\begin{alignat*}{2}
\frac{\partial a_{11}}{\partial \theta} & = -a_{21} & \qquad 
\frac{\partial a_{12}}{\partial \theta} & = -a_{22} \\[0.1in]
\frac{\partial a_{21}}{\partial \theta} & = a_{11} - I a_{21} & 
\qquad \qquad \frac{\partial a_{22}}{\partial \theta} & = a_{12} - I 
a_{22}.
\end{alignat*}
In other words, the function pairs $(a_{11}, a_{21})$ and $(a_{12}, 
a_{22})$ are each solutions of the system
\begin{align}
\frac{\partial f}{\partial \theta} & = -g \label{ODEsys} \\[0.1in]
\frac{\partial g}{\partial \theta} & = f - I g \notag
\end{align}
for functions $f(x,y,\theta),\, g(x,y,\theta)$.
The solution of these differential equations depends on the value of $I$.

\subsubsection{Case 1: $I^2 > 4$}

The general solution of \eqref{ODEsys} in this case is
\begin{align*}
f & = c_1(x,y) e^{r_1\theta} + c_2(x,y) e^{r_2 \theta} \\
g & = -c_1(x,y) r_1 e^{r_1\theta} - c_2(x,y) r_2 e^{r_2 \theta},
\end{align*}
where
\[ r_1, r_2 = \tfrac{1}{2}(-I \pm \sqrt{I^2 - 4}). \]
By modifying $x$ and $y$ if necessary, we can assume that
\begin{alignat*}{2}
a_{11} & = c_1(x,y) e^{r_1\theta} & \qquad a_{12} & = c_2(x,y) e^{r_2\theta} \\
a_{21} & = -c_1(x,y) r_1 e^{r_1\theta} & \qquad a_{22} & = -c_2(x,y) 
r_2 e^{r_2\theta}.
\end{alignat*}
Then the first two equations of \eqref{homogstruct} imply that $c_1$ 
is a function of $x$ alone and $c_2$ is a function of $y$ alone.

Now the third equation of \eqref{homogstruct} implies that
\[ d\eta^3 \equiv c_1(x) c_2(y) e^{-I\theta} \sqrt{I^2-4}\, dx \& dy 
\mod{\eta^3}. \]
By the Pfaff theorem, there exists a function $z$ on $\Sigma$, 
independent from $x, y, $ and $\theta$, such that
\[ \eta^3 = c_1(x) c_2(y) e^{-I\theta} \sqrt{I^2-4} \left( dz + 
\tfrac{1}{2}(x\, dy - y\, dx) \right). \]
Finally, the third equation of \eqref{homogstruct} now implies that 
$c_1$ and $c_2$ are in fact constants.  Without loss of generality, 
we may assume that $c_1 = c_2 = 1$, and that our coframing has the 
form
\begin{align*}
\eta^1 & = e^{r_1\theta}\, dx + e^{r_2\theta}\, dy \\
\eta^2 & = -r_1 e^{r_1\theta}\, dx - r_2 e^{r_2\theta}\, dy \\
\eta^3 & = e^{-I\theta} \sqrt{I^2-4} \left( dz + \tfrac{1}{2}(x\, dy 
- y\, dx) \right) \\
\phi & = d\theta.
\end{align*}

\subsubsection{Case 2: $I^2 < 4$}

The general solution of \eqref{ODEsys} in this case is
\begin{align*}
f & = e^{-I\theta/2} [c_1(x,y) \cos(r\theta) + c_2(x,y) \sin(r\theta) ] \\
g & = \tfrac{1}{2}e^{-I\theta/2} [c_1(x,y) (I \cos(r\theta) + 
r\sin(r\theta)) \\
& \qquad \qquad \qquad + c_2(x,y) (I \sin(r\theta) - r \cos(r\theta)) ],
\end{align*}
where
\[ r = \tfrac{1}{2}\sqrt{4 - I^2}. \]
A similar argument to that given above shows that that we can take 
our coframing to be
\begin{align*}
\eta^1 & = e^{-I\theta/2} [\cos(r\theta)\, dx + \sin(r\theta)\, dy] \\
\eta^2 & = \tfrac{1}{2}e^{-I\theta/2} [(I \cos(r\theta) + 
r\sin(r\theta))\, dx + (I \sin(r\theta) - r \cos(r\theta))\, dy ] \\
\eta^3 & = -r e^{-I\theta}  \left( dz + \tfrac{1}{2}(x\, dy - y\, dx) 
\right) \\
\phi & = d\theta.
\end{align*}

\subsubsection{Case 3: $I = 2$}
The general solution of \eqref{ODEsys} in this case is
\begin{align*}
f & = e^{-t} [-c_1(x,y)\, t + c_2(x,y)(1 + t)] \\
g & = e^{-t} [c_1(x,y)(1 - t) + c_2(x,y)\, t].
\end{align*}
A similar argument to that given above shows that that we can take 
our coframing to be
\begin{align*}
\eta^1 & = e^{-\theta} [(1 + \theta)\, dx - \theta\, dy] \\
\eta^2 & = e^{-\theta} [\theta\, dx + (1 - \theta)\, dy ] \\
\eta^3 & = e^{-2\theta}  \left( dz + \tfrac{1}{2}(x\, dy - y\, dx) \right) \\
\phi & = d\theta.
\end{align*}

\subsubsection{Case 4: $I = -2$}
The general solution of \eqref{ODEsys} in this case is
\begin{align*}
f & = e^{t} [-c_1(x,y)\, t + c_2(x,y)(1 - t)] \\
g & = e^{t} [c_1(x,y)(1 + t) + c_2(x,y)\, t].
\end{align*}
A similar argument to that given above shows that that we can take 
our coframing to be
\begin{align*}
\eta^1 & = e^{\theta} [(1 - \theta)\, dx - \theta\, dy] \\
\eta^2 & = e^{\theta} [\theta\, dx + (1 + \theta)\, dy ] \\
\eta^3 & = e^{2\theta}  \left( dz + \tfrac{1}{2}(x\, dy - y\, dx) \right) \\
\phi & = d\theta.
\end{align*}

Note that none of these four coframings have coordinate expressions 
which are periodic in the $\theta$ variable.  Consequently, in all 
four cases the indicatrix for the sub-Finsler metric in each tangent 
space fails to be a closed curve.  (In fact, these indicatrices are 
not even connected.)  Therefore, these sub-Finsler structures exist 
``micro-locally" -- that is, in a neighborhood of each point in 
$\Sigma$ -- but not locally.  In other words, there is no open set $\scU 
\subset \scX$ for which the sub-Finsler metric is defined on all of 
$T\scU$.  This is consistent with a theorem of Rund (see 
\cite{BCS00}) which implies that for any Minkowski norm on a plane 
$D_x$, the average value of $I$ over the indicatrix must be zero; 
therefore, if $I$ is any nonzero constant, the indicatrix cannot 
possibly be a closed, strongly convex curve in $D_x$.

\subsection{Symmetry groups of dimension three}

Now we consider the more inclusive case where the group $G$ of
symmetries of $\Sigma$ is three-dimensional and acts transitively on $\scX$.
Since $\Sigma$ is invariant under the action of $G$, it is completely
determined by the fiber $\Sigma_x$ at any point $x \in \scX$.
Conversely, if we fix a point $x \in \scX$ and choose any smooth curve
$\Gamma \subset D_x$ which is strongly convex and encloses the origin,
then there exists a unique sub-Finsler structure $\Sigma$ on $(\scX, D)$
which is invariant under the action of $G$ and satisfies $\Sigma_x =
\Gamma$.

Thus we conclude that the sub-Finsler structures of this
type are generated by choosing a three-dimensional Lie group $G$ (or 
a quotient thereof by a discrete subgroup), a
2-plane $D \subset T_eG$ which is not a Lie subalgebra (so that it is
bracket-generating), and a smooth curve $\Gamma$ in $D$ which is
strongly convex and surrounds the origin.

\begin{example}\label{Heisenberg}
Let $\scH$ be the Heisenberg group, defined by
\begin{equation*}
   \scH = \left\{ \begin{bmatrix}
                1 & y & z + \frac{1}{2}xy \\
                0 & 1 & x \\
                0 & 0 & 1
             \end{bmatrix} : x,y,z \in \bb{R}  \right\} \cong \bb{R}^{3},
\end{equation*}
and let the contact structure on $\scH$ be the rank two distribution
\begin{equation}\label{pfaffnormalform}
    D = \{dz + \frac{1}{2}(x\,dy - y\,dx)\}^{\perp}.
\end{equation}

The existence of this global coordinate system makes 
it easy to describe sub-Finsler geodesics within the Heisenberg group. 
Moreover, this example is prototypical:  by a theorem of Pfaff 
(see \cite{BCG3}), any contact 3-manifold has local coordinates 
$(x,y,z)$ for which the contact system is given by the symmetric 
normal form \eqref{pfaffnormalform} above.  

We can define a homogeneous, flat sub-Riemannian metric on $(\scH, 
D)$ by declaring the vectors
\[ V_1 = \frac{\partial}{\partial x} + \frac{y}{2} 
\frac{\partial}{\partial z}, \qquad V_2 = \frac{\partial}{\partial y} 
- \frac{x}{2} \frac{\partial}{\partial z} \]
to be orthonormal.  Let $\Sigma_1$ be the unit circle bundle for this 
sub-Riemannian structure on $\scH$, and define a coordinate $\theta$ 
on $\Sigma_1$ by the condition that, for $u \in \Sigma_1$,
\[ u = (\cos \theta) V_1 + (\sin \theta) V_2. \]

It is straightforward to check that $V_1, V_2$ are left-invariant, 
horizontal vector fields on $(\scH, D)$, and that therefore any 
scaling function $r(\theta)$ which depends on $\theta$ alone defines 
a homogeneous sub-Finsler structure on $\scH$.  It is also 
straightforward to check that the coframing
\begin{align*}
\bar{\eta}^1 & =\rho^*\left( (r \cos \theta - r'\sin \theta)\,dx - (r 
\sin \theta + r'\cos \theta)\,dy\right) \\
\bar{\eta}^2 & = \rho^*\left( \sqrt{r(r + r'')}[(\sin \theta)\, dx + 
(\cos \theta)\, dy] \right)\\
\bar{\eta}^3 & = \rho^*\left( r^{3/2}\sqrt{r + r''}[dz + 
\tfrac{1}{2}(x\,dy - y\,dx)] \right)\\
\bar{\phi} & = \rho^*\left( \frac{\sqrt{r + r''}}{\sqrt{r}}\, d\theta\right)
\end{align*}
on the sub-Finsler structure $\Sigma$ defined by $r(\theta)$ is 
4-adapted.  The invariants for this coframing are
\begin{gather*}
I = -\frac{1}{2}\frac{(rr'''  + 3r'r'' + 4rr')}{\sqrt{r}(r + r'')^{3/2}}, \\
K = A_1 = A_2 = J_1 = J_2 = S_0 = S_1 = S_2 = 0.
\end{gather*}
In this case, the geodesic equations \eqref{geodesic_eqns} can be written as
\begin{align}
dx & = \frac{\cos\theta(s)}{r(\theta(s))}\, ds \notag \\[0.1in]
dy & = -\frac{\sin\theta(s)}{r(\theta(s))}\, ds \notag \\[0.1in]
dz & = \frac{x(s)\sin\theta(s) +  y(s)\cos\theta(s)}{2r(\theta(s))}\, 
ds \label{Heisenberg_geo_eqs} \\[0.1in]
d\theta & = \frac{\sqrt{r(\theta(s))}}{\sqrt{r(\theta(s)) + 
r''(\theta(s))}} \, \lambda(s)\, ds \notag \\[0.1in]
d\lambda & = I \lambda^2(s)\, ds. \notag
\end{align}
Since the Lie algebra of $\scH$ is solvable, it is no 
surprise to find that these equations can be solved by quadrature. 
First we prove a straightforward but useful 
lemma.

\begin{lemma}\label{positivity}  The expression $r(r + r'')$ is positive and bounded away from 
zero.
\end{lemma}
\begin{proof}
Let $x$ be any point in $\scH$, and 
let $\Sigma_{x} \subset D_{x}$ be the indicatrix in $D_{x}$ 
corresponding to a sub-Finsler metric $F$.  The strong convexity of 
$F$ is equivalent to the condition that, for any parametrization 
$(u(t), v(t))$ of 
$\Sigma_{x}$,
\begin{equation}\label{strongconvexity}
  \frac{u'' 
v' - u' v''}{u' v - u v'} > 0 
\end{equation}
everywhere on $\Sigma_x$
(see \cite{BCS00} for a 
proof).  In particular, for the parametrization 
\[ u(\theta) = R(\theta)
\cos{\theta}, \qquad v(\theta) = R(\theta) \sin{\theta}, \]
this is equivalent to 
saying that
\begin{equation}\label{strongconvexity2}
\frac{r + r''}{r}  =  \frac{R^{2} + 2(R')^{2} - RR''}{R^{2}} > 0 
\end{equation}
everywhere on $\Sigma_x$, and hence $r(r + r'') > 0$ as well.
(Recall that $r$ is the 
\emph{reciprocal} of the radial position function $R$ defining $\Sigma$.) 
Since $\Sigma_{x}$ is compact, this quantity is bounded away from 
zero.  By the homogeneity of $\scH$, this bound is the same at every 
point $x \in \scH$.
\end{proof}

This lemma and the geodesic 
equations \eqref{Heisenberg_geo_eqs} imply the following result.
\begin{theorem}
  For any 
homogeneous sub-Finsler metric $F$ on the Heisenberg group $\scH$, 
the sub-Finsler geodesics are straight lines parallel to the 
$xy$-plane or liftings of simple closed curves in the $xy$-plane.  In 
the latter case, the simple closed curves are the curves of minimal 
Finsler arc length enclosing a given Euclidean area in the 
plane.
\end{theorem}
\begin{proof}
A curve $\gamma:[a,b] \to \scH$ is a geodesic if and only if it is an integral curve of the system \eqref {Heisenberg_geo_eqs}.
The equations for $d\theta$ and 
$d\lambda$ in  \eqref{Heisenberg_geo_eqs} allow us to 
write
\begin{align}
  \frac{d\lambda}{\lambda} &= I \sqrt{\frac{r + 
r''}{r}}\, d\theta \notag\\
  &= -\frac{1}{2}\frac{r r''' + 3r' r'' + 
4r r'}{r(r + r'')} \, d\theta \notag\\
  &= -\frac{1}{2}\frac{d(r(r + 
r''))}{r(r + r'')} - \frac{dr}{r}, \notag
\end{align}
and 
so
\begin{equation*}
  \lambda = \frac{c \lambda_{0}}{r \sqrt{r (r + 
r'')}},
\end{equation*}
where $\lambda_0,\, c$ are constants, with $c 
= r(0)\sqrt{r(0)(r(0) + r''(0))} > 0$.

Integral curves of \eqref{Heisenberg_geo_eqs} corresponding to
$\lambda_{0} = 0$ are 
straight lines parallel to the $xy$-plane.  If $\gamma$ is an integral curve corresponding to some $\lambda_{0} \neq 0$, 
then we have
\begin{equation*}
  d\theta = \frac{c \lambda_{0}}{r(r + r'')} 
ds.
\end{equation*}
By the preceding lemma, the quantity in the denominator is positive 
and bounded away from zero; thus $\theta$ 
varies monotonically with $s$, without bound.  We may therefore 
reparametrize the equations for $dx$, $dy$ and $dz$ in terms of 
$\theta$.  If $(x_{0},\, y_{0}) = (x(\theta_0),\, y(\theta_0))$ is 
any point on the projection of $\gamma$ to the $xy$-plane, 
then for any other value $\theta$, we have
\begin{align}
  x(\theta) 
- x_{0} &= \frac{1}{c \lambda_{0}} \int_{\theta_{0}}^{\theta} (r + 
r'') \cos{t} \, dt \label{xyeqs} \\
  y(\theta) - y_{0} &= 
-\frac{1}{c \lambda_{0}} \int_{\theta_{0}}^{\theta} (r + r'') \sin{t} 
\, dt .\notag
\end{align}
Integrating by parts twice shows that 
\begin{align}
x(\theta) - x_{0} &= \frac{-1}{c \lambda_{0}} \left(r(\theta)^{2} u'(\theta) - r(\theta_{0})^{2} u'(\theta_{0})\right) \notag\\
y(\theta) - y_{0} &= \frac{1}{c \lambda_{0}} \left(r(\theta)^{2} v'(\theta) - r(\theta_{0})^{2} v'(\theta_{0})\right) ,\notag
\end{align}
where 
\[ u(\theta) = R(\theta)\cos{\theta}, \qquad v(\theta) = R(\theta)\sin{\theta} \] 
is the parametrization of the indicatrix used in Lemma \ref{positivity}.  Thus $x(\theta) - x_{0}$ and $y(\theta) - y_{0}$ are simultaneously zero if and only if
\begin{equation*}
  \begin{bmatrix}
      u'(\theta)\\[0.1in]
      v'(\theta)
   \end{bmatrix} = \frac{r(\theta_{0})^{2}}{r(\theta)^{2}} \begin{bmatrix}
                                                                                                        u'(\theta_{0})\\[0.1in]
                                                                                                        v'(\theta_{0})
                                                                                                  \end{bmatrix},
\end{equation*}
and because the indicatrix is strongly convex, this occurs precisely when $\theta = \theta_{0} + 2n\pi$ for any integer $n$.  Since $\theta$ is unbounded as a function of arc length, it attains these values.  Therefore, the projection of $\gamma$ onto the $xy$-plane is a simple closed curve.

Finally, since $dz = - \frac{1}{2}(x\,dy - y\,dx)$, Stokes' theorem implies that the difference $z(\theta) - z_{0}$ along any 
horizontal curve in $\scH$ is proportional to the signed area 
enclosed by the projection of the curve onto the $xy$-plane and the 
line segment connecting $(x(\theta),y(\theta))$ to $(x_{0}, y_{0})$. 
Thus $z$ varies monotonically with increasing $\theta$, and the 
projection of $\gamma$ onto the $xy$-plane is the curve of 
shortest Finsler arc length enclosing a given Euclidean area in the 
plane.
\end{proof}

By the last observation in the proof above, the 
projections of sub-Finsler geodesics are precisely the solutions to 
the dual of the classical isoperimetric problem known as \emph{Dido's 
problem}, named after Queen Dido in Virgil's \emph{Aeneid} 
(see \cite{Montgomery02}).  The classical solutions (using Riemannian arc 
length) are circles.  In the Finsler case, the solution curves need 
not be circles, as one of the examples below 
illustrates.
\end{example}

\subsubsection{Randers metrics on $\scH$ and their geodesics
}  Consider the Randers-type, homogeneous sub-Finsler metric 
on $D$ obtained by choosing the function $r(\theta)$ to be
\begin{equation*}
   r(\theta) = 1 + B \cos{\theta}, \,\, 0 < B < 
1.
\end{equation*}
The indicatrix $\Sigma_x$ for this metric in each plane $D_x$ is the 
off-center ellipse with polar equation
\[ R = \frac{1}{1 + B \cos{\theta}}. \]
For this metric,
\[ I = \frac{3 B \sin\theta}{2\sqrt{1 + B\cos\theta}}.\] 

The geodesic equations \eqref{Heisenberg_geo_eqs} can be integrated 
explicitly in terms of $\theta$.  When $\lambda_{0} = 0$, the 
geodesics are straight lines; when $\lambda_0 \neq 0$, integrating 
yields:
\begin{align}
x(\theta) - x_{0} & = \frac{1}{k}(\sin{\theta} - \sin{\theta_{0}})\notag \\[0.1in]
y(\theta) - y_{0} & = \frac{1}{k}(\cos{\theta} - 
\cos{\theta_{0}})\label{Randers_geo_eqs} \\[0.1in]
z(\theta) - z_{0} & = \frac{1}{2k^{2}}[\theta - \theta_{0} - \sin{(\theta - 
\theta_{0})}]  + \frac{1}{2k}[y_{0}(\sin{\theta} - \sin{\theta_{0}}) 
- x_{0}(\cos{\theta} - \cos{\theta_{0}})]\notag,
\end{align}
where $k 
= \lambda_{0}\,\sqrt{(1 + B\cos{\theta_{0}})^{3}}$.  These geodesics 
are liftings of circles of radius $1/k$ in the $xy$-plane.  In the 
limiting case $B = 0$, we recover the geodesics of the flat 
sub-Riemannian metric.  When $0 < B < 1$, the anisotropy of the 
indicatrix is manifested in the way that the area of the projected 
circle in the $xy$-plane (and, therefore, $dz/d\theta$) varies with 
the initial value $\theta_{0}$, unlike in the sub-Riemannian 
case.

Figure 1 shows some typical geodesics for $B = 
\frac{1}{2}$, starting from $(x_{0},y_{0},z_{0})=(0,0,0)$, with 
initial values $\lambda_0 = 0.3$ and $\theta_0 = 0, \pi/2$, and 
$\pi$.  For comparison, Figure 2 shows geodesics for the flat 
sub-Riemannian metric with these same initial values.

\noindent
\begin{center}
\begin{minipage}{.4\linewidth}
\centering\epsfig{figure=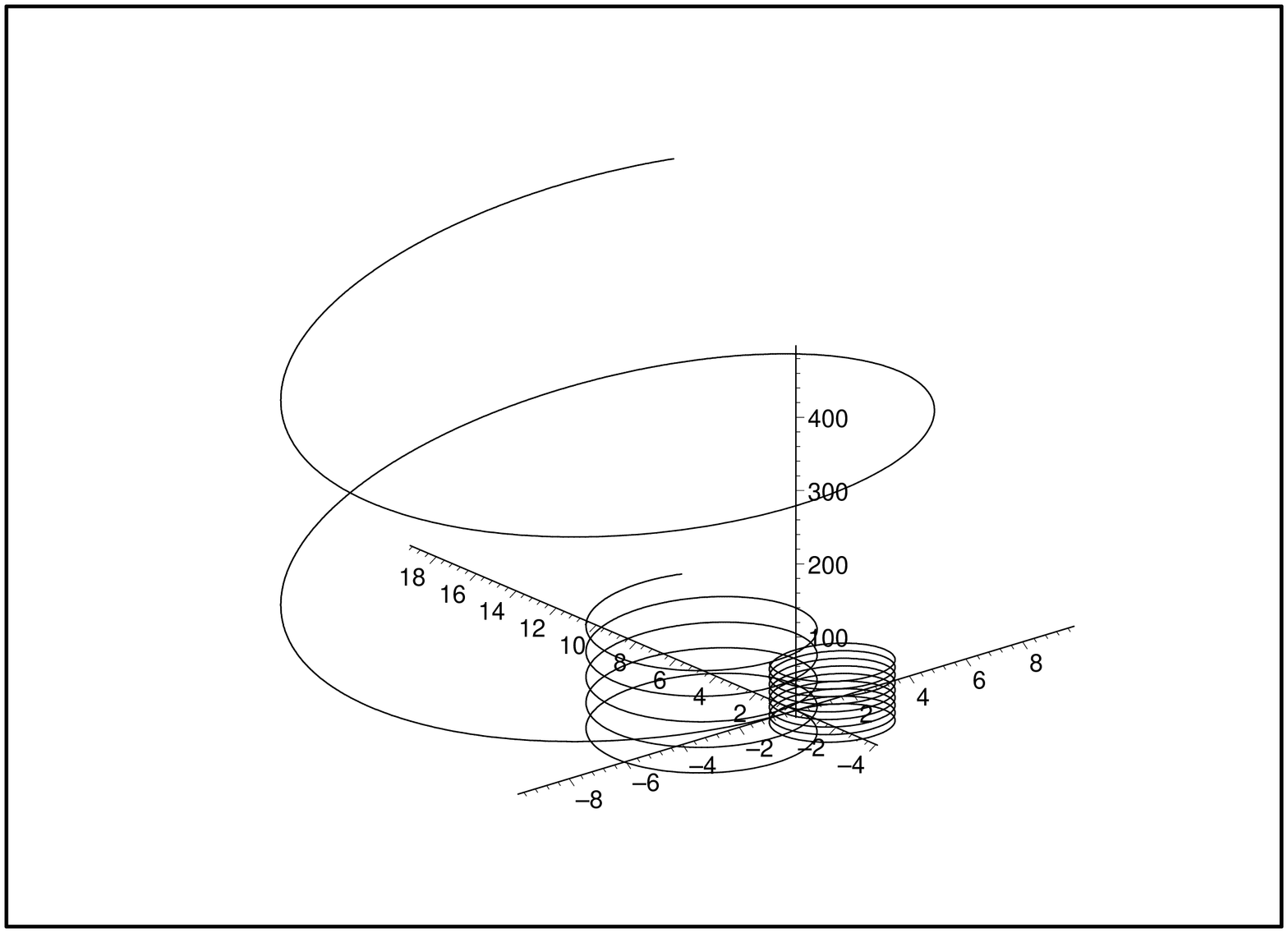,width=\linewidth,angle=-90}
\newline\caption{Figure 1: Geodesics of the Randers metric}
\end{minipage}\hspace{.6in}
\begin{minipage}{.4\linewidth}
\centering\epsfig{figure=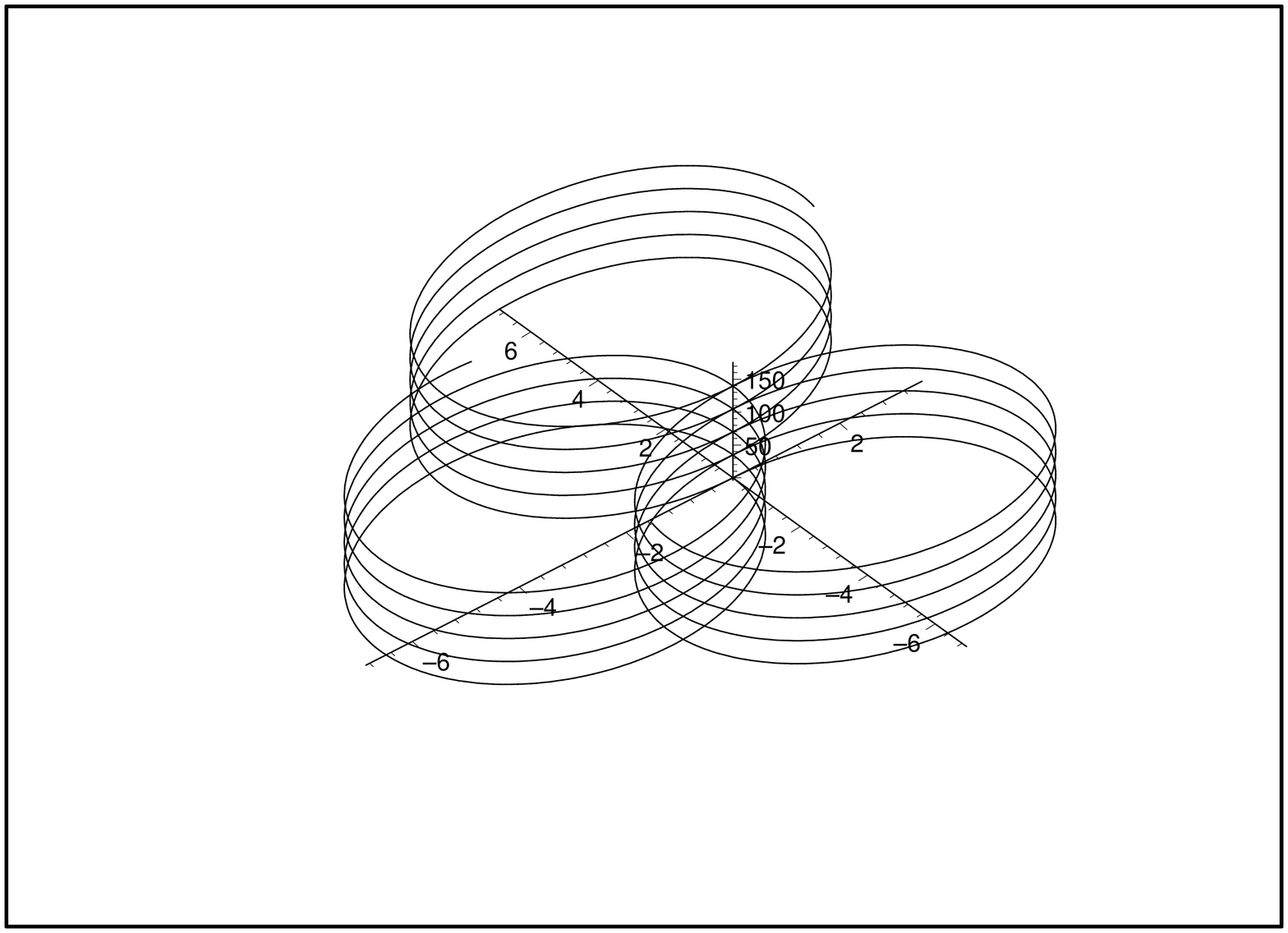,width=\linewidth,angle=-90}
\newline\caption{Figure 2: Geodesics of the sub-Riemannian metric}
\end{minipage}
\end{center}

\vspace{0.2in}

\subsubsection{A
``lima\c{c}on metric" on $\scH$ and its geodesics}
For an example in which the 
geodesics are not liftings of conic sections, consider the 
sub-Finsler metric whose indicatrix is the convex lima\c{c}on with 
polar equation $R = 3 + \cos{\theta}$,
so that
\[r(\theta) = 
\frac{1}{3 + \cos{\theta}}. \]
For this metric,
\[I = \frac{-3 
\sin{\theta}(15 \cos{\theta} + 13)}{2\sqrt{(9\cos{\theta} + 
11)^{3}}}.\]
As always, the geodesics are straight lines when $\lambda_{0} = 0$, but otherwise they are liftings of curves in the 
$xy$-plane defined by the equations
\begin{align}
x(\theta) - x_{0} & = \frac{1}{2L} \left( \frac{\sin{\theta}(4\cos{\theta} + 
6)} {(3 + \cos{\theta})^{2}} - 
\frac{\sin{\theta_{0}}(4\cos{\theta_{0}} + 6)} {(3 + 
\cos{\theta_{0}})^{2}} \right)\notag  \\[0.1in]
y(\theta) - y_{0} & = - \frac{1}{L} \left( \frac{9\cos{\theta} + 19}{(3 + 
\cos{\theta})^{2}} -  \frac{9\cos{\theta_{0}} + 19}{(3 + 
\cos{\theta_{0}})^{2}} \right), \notag \end{align}
where
\[L = 
\lambda_{0}\frac{\sqrt{9\cos{\theta_{0}} + 11}}{(3 + 
\cos{\theta_{0}})^{3}}, \,\, \lambda_{0} \neq 0.\]

These curves in 
the $xy$-plane are not circles, nor are they ellipses (or 
lima\c{c}ons).  Figure 3 shows geodesics for this metric starting 
from $(x_{0},y_{0},z_{0})=(0,0,0)$, with initial values $\lambda_{0} 
= 1$ and $\theta_0 = 0, \pi/2$, and $\pi$.   Figure 4 shows the 
projections of these curves onto the 
$xy$-plane.

\noindent
\begin{center}
\begin{minipage}{.4\linewidth}
\centering\epsfig{figure=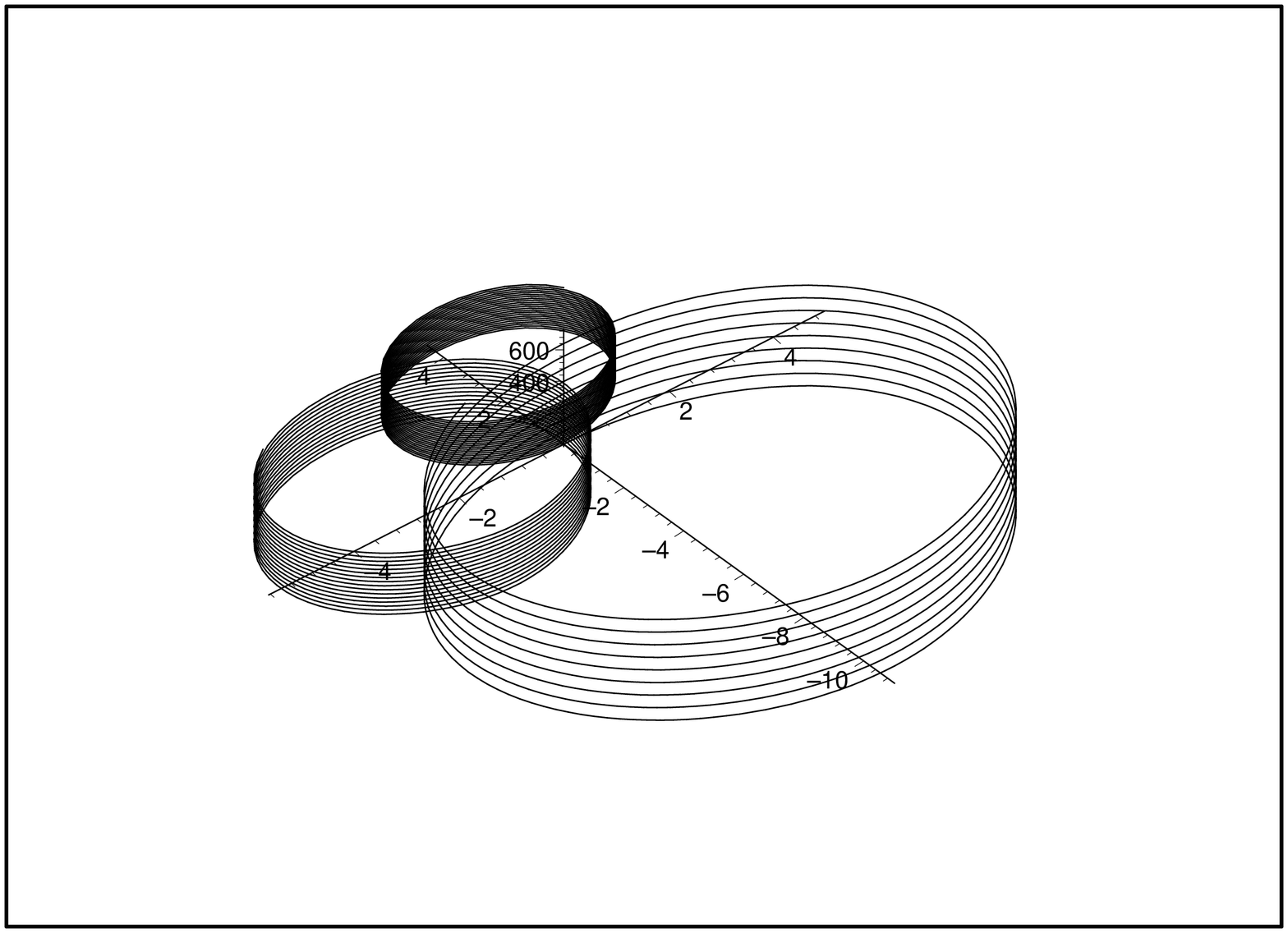,width=\linewidth,angle=-90}
\newline\caption{Figure 3: Geodesics of the lima\c{c}on metric}
\end{minipage}\hspace{.6in}
\begin{minipage}{.4\linewidth}
\centering\epsfig{figure=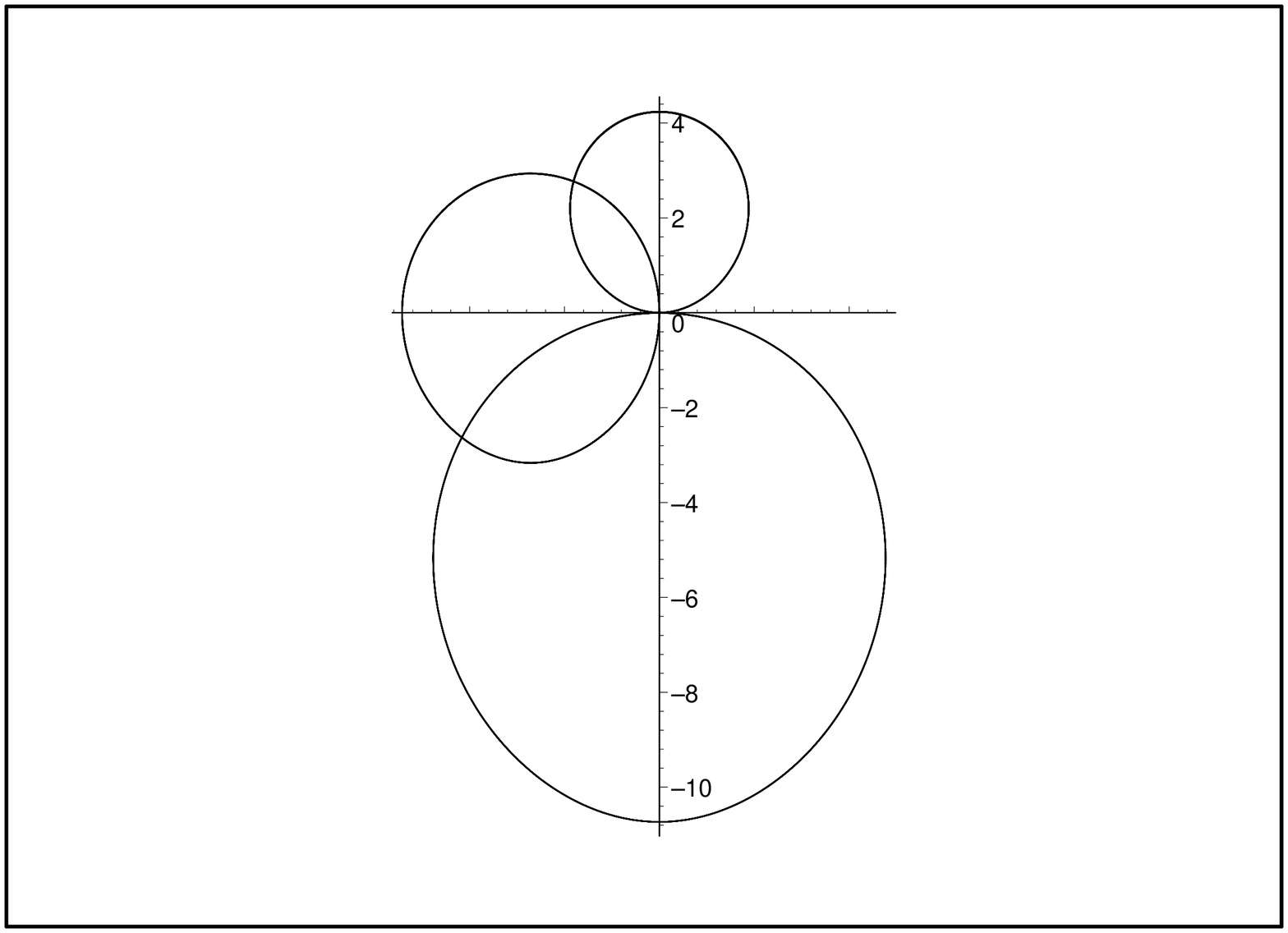,width=\linewidth,angle=-90}
\newline\caption{Figure 4:  Projections of lima\c{c}on metric 
geodesics onto the $xy$-plane}
\end{minipage}
\end{center}

\vspace{0.2 in}

\section{Conclusion}

We have only begun to explore sub-Finsler geometry in this paper, and 
we have every reason to believe that it will become a useful 
extension of sub-Riemannian geometry, particularly in the context of 
control theory.  In future papers, we plan to investigate 
higher-dimensional cases (including the important phenomenon of 
abnormal geodesics), singularities, and other topics likely to be of 
interest for control theory applications.

\end{document}